%% file: GCRODR-direct-proj.tex
\documentclass[final]{siamltex}
\usepackage{euscript,graphicx,epsfig,amsfonts,mathtools,epstopdf}
\usepackage[ruled,lined,vlined,linesnumbered,algosection]{algorithm2e}
\usepackage{amssymb}
\usepackage{enumerate}
\usepackage[dvips,
          letterpaper=true,
          colorlinks=true,
          linkcolor=red,
          filecolor=green,
          citecolor=red,
          pdfpagemode=None]
          {hyperref}

\input{macros}
\def\RGMRES{Recycled GMRES}
\newcommand*{\Resize}[2]{\resizebox{#1}{!}{$#2$}}%
\newcommand*{\Sf}[1]{\Resize{.03\textwidth}{\left(#1\right)}}
\def\Sfs{\Sf{\sigma}}
\newcommand*{\SfAlg}[1]{\Resize{.04\textwidth}{\left(#1\right)}}

\title{Two recursive GMRES-type methods for shifted linear systems with general preconditioning\thanks{%
This version dated \today.}}

\author{Kirk M. Soodhalter\footnotemark[2]}

\begin{document}
\maketitle
\renewcommand{\thefootnote}{\fnsymbol{footnote}}�
\footnotetext[2]{Industrial Mathematics Institute, Johannes Kepler University, 
Altenbergerstra{\ss}e 69, A-4040 Linz, Austria.
({\tt kirk.soodhalter@indmath.uni-linz.ac.at})}
 
\begin{abstract}
We present two minimum residual methods for solving sequences of 
shifted linear systems, 
the right-preconditioned shifted GMRES and shifted
\RGMRES\  algorithms which use a seed projection strategy often employed to solve multiple
related problems.
These methods
are compatible with general preconditioning of all systems, and when
restricted to right preconditioning, require no extra applications of the 
operator or preconditioner.  These seed projection methods perform a minimum
residual iteration for the base system while improving the approximations
for the shifted systems at little additional cost.
The iteration continues until the base system approximation
is of satisfactory quality.  The method is then recursively called for the
remaining unconverged systems.  
We present both methods inside of a general framework which allows these
techniques to be extended to the setting of flexible preconditioning and 
inexact Krylov methods.
We present some analysis 
of such methods and numerical experiments demonstrating the effectiveness
of the algorithms we have derived.
\end{abstract}



\pagestyle{myheadings}
\thispagestyle{plain}
\markboth{K. M. SOODHALTER}{RIGHT PRECONDITIONED METHODS FOR NONHERMITIAN SHIFTED SYSTEMS}

\section{Introduction}\label{section.intro}
We develop techniques for solving a family (or a sequence of families) of linear systems
in which the coefficient matrices differ only by a scalar multiple of the identity.  
There are many applications which warrant the solution 
of a family of shifted linear 
systems,  such as  
those arising in lattice quantum chromodynamics (QCD) 
(see, e.g., \cite{Frommer1995}) as well as other applications such as 
Tikhonov-Phillips regularization, global methods of nonlinear analysis, and Newton trust region 
methods \cite{CDT1985}.
The goal is to develop a framework in which minimum residual
methods can be applied to shifted systems in a way that:
\begin{enumerate}[(a)]
\item\label{item.shift-rel} allows us to exploit the relationships 
between the coefficient matrices
\item\label{item.gen-precond} is compatible with general (right) preconditioning. 
\end{enumerate}

In this paper, we use such a framework to propose two new methods: one which is built on top
of the GMRES method \cite{Saad.GMRES.1986} for solving a family of shifted systems (cf. \eqref{eqn.one-shifted-fam})
and one which is built on top of a GCRO-type augmented Krylov method \cite{deSturler.GCRO.1996} which,
when paired with a harmonic Ritz vector recycling strategy \cite{Morgan1998,Paige1995}, 
is an extension of the GCRO-DR method \cite{Parks.deSturler.GCRODR.2005} to solve a sequence of shifted system families 
(cf. \eqref{eqn.one-shifted-family-seq}).  To do this, we use a seed projection strategy, often proposed for use
in conjunction with short-term recurrence iterative methods 
\cite{CN.GalProjAnalMRHS.1999,Chan1997,KMR.2001,Saad1987}.

The rest of this paper is organized as follows.  In the next section, we discuss some previous strategies to
treat such problems and discuss some of their limitations.
In Section \ref{section.prelim}, we review
the minimum residual Krylov subspace method GMRES as well as two GMRES 
variants, one for shifted linear systems and the other extending GMRES to
the augmented Krylov subspace setting, i.e., \RGMRES.  In Section \ref{section.dir-proj-framework}, 
we present a general framework to perform minimum residual projections
of the shifted system residuals with respect to the search space generated
for the base system.  In Subsection \ref{subsection.GMRES-shift} we use this framework to derive our
\textbf{shifted GMRES} method and in Subsection \ref{subsection.rgmres-shift} we derive a \textbf{shifted \RGMRES}\  method.
In Section \ref{section.analysis}, we present some analysis of the expected performance of these
methods.  In Section \ref{section.num-results}, we present
some numerical results before concluding in Section \ref{section.conclusions}.

\section{Background}\label{section.background}
Consider a family of shifted linear systems, which we parameterize
by $\ell$, i.e.,  
\be\label{eqn.one-shifted-fam}
	\prn{\vek A+\sigma^{(\ell)}\vek I}\vek x^{(\ell)}=\vek b\mbox{\ \ \ for\ \ \ }\ell=1,\ldots , L.
\ee
We call the numbers $\curl{\sigma^{(\ell)}}_{\ell=1}^{L}\subset\C$ \textit{shifts}, $\vek A$ the \textit{base matrix}, and $\vek A+\sigma\vek I$ a \textit{shifted matrix}.  
Systems of the form (\ref{eqn.one-shifted-fam}) are called \textit{shifted linear systems}.  
 Krylov subspace methods have been 
proposed to simultaneously solve this family of systems, see, e.g., \cite{Darnell2008,Frommer2003, Frommer1998,Kilmer.deSturler.tomography.2006,Simoncini2003a}.  These methods satisfy requirement
(\ref{item.shift-rel}) but are not compatible with general preconditioning strategies,
as they rely on the invariance of the Krylov subspace under constant 
shift of the coefficient matrix; cf.  (\ref{eqn.shift-invariance}).
Specially chosen polynomial preconditioners, however, have been shown to be compatible
with such methods; see, e.g., \cite{ASV.2012,BKLSS.mpgmres.2016,BG.2014,Jegerlehner.shifted-systems.1996,M.2013,SBK.2014, WWJ.2012}.

 We can introduce an additional parameter $i$, which indexes a sequence of matrices 
$\curl{\vek A_{i} }\subset\Cnn$, and for each $i$, we solve a family of the form
\be\label{eqn.one-shifted-family-seq}
	\prn{\vek A_{i}+\sigma_{i}^{(\ell)}\vek I}\vek x_{i}^{(\ell)} = \vek b_{i}\mbox{\ \ \ for\ \ \ }\ell=1\ldots L_{i}
\ee
We consider the case that the right-hand side varies with respect to $\vek A_{i}$ but not for
each shift.  What we propose is indeed applicable in the more general setting, but we do not treat that here.
Augmented Krylov subspace methods have been proposed for 
efficiently solving a sequence of linear systems with 
a slowly changing coefficient matrix, allowing important spectral information 
generated while solving $\vek A_{i}\vek x_{i} = \vek b_{i}$ to be used to augment the
Krylov subspace generated when solving
$\vek A_{i+1}\vek x_{i+1} = \vek b_{i+1}$; see, e.g.,  
\cite{Parks.deSturler.GCRODR.2005,SYEG.2000,WSG.2007}.  In cases such as a Newton
iteration, these matrices are available one at a time, while in a case such as
an implicit time-stepping scheme, the matrix may not change at all.

In \cite{SSX.2014}, the authors explored solving a family of shifted systems over an 
augmented Krylov subspace.  Specifically, the goal was to develop a method which solved 
the family of systems simultaneously, using one augmented 
subspace to extract all candidate solutions, 
which also had a fixed storage requirement, independent of the number of shifts $L$.  
It was shown
that in general within the framework of GMRES for shifted systems \cite{Frommer1998} 
and subspace recycling
\cite{Parks.deSturler.GCRODR.2005}, such a method, \textit{does not exist}. 
In the context of subspace recycling for Hermitian linear systems, in the absence
of preconditioning Kilmer and de Sturler proposed a MINRES method in a subspace recycling 
framework which 
simultaneously solves multiple non-Hermitian systems, which all differ from a real-symmetric system by a
complex multiple of the identity \cite{Kilmer.deSturler.tomography.2006}, by minimizing the shifted residuals over the augmented Krylov subspace subspace, built using the symmetric Lanczos process.  
In this paper, we 
focus exclusively on problems in which the base coefficient matrices $\vek A_{i}$ are
\emph{non-Hermitian}.

A conclusion one can draw from \cite{SSX.2014} is that we should consider avoiding methods
relying on the invariance of Krylov subspaces under a constant shift of the identity; cf.\ \eqref{eqn.shift-invariance}.  Relying on this invariance imposes restrictions on 
our ability to develop an algorithm.  Furthermore, relying on this shift invariance means we cannot 
use arbitrary preconditioners.  
General preconditioners are unavailable if we want to exploit shift invariance, as Krylov subspaces generated by preconditioned
systems are not invariant with respect to a shift in the coefficient matrix. 
In the case that preconditioning is not used, a subspace recycling technique 
has been proposed \cite{S.2016}, built on top of the Sylvester equation interpretation of
\eqref{eqn.one-shifted-fam} observed by Simoncini in \cite{S.1996}.  However, this is also not compatible
with general preconditioning.

Learning from the results in \cite{SSX.2014}, we focus on methods which do
not rely on the shift invariance.  Rather than focusing on specific Krylov subspace techniques
(augmented or not), we instead begin by developing 
a general framework of minimum residual
projection techniques for shifted linear systems.     
In this framework,
we extract candidate solutions for all shifted systems from the augmented Krylov subspace for one linear system and we 
select each candidate solution according to a minimum residual Petrov-Galerkin condition.  
This framework is compatible with arbitrary right preconditioners, and the computational 
cost for each
additional shifted system is relatively small but nontrivial.  
By specifying subspaces once the framework is developed, we derive
minimum residual
methods for shifted systems that are compatible
with general right preconditioning.  Though not considered in this paper, 
the framework is also compatible with flexible and inexact Krylov methods.
These methods descend from the
Lanczos-Galerkin seed methods, see, e.g., \cite{CN.GalProjAnalMRHS.1999,Chan1997,KMR.2001,Saad1987}.

In this work, we restrict ourselves to right preconditioned
methods.  Doing this allows us to derive methods
which require extra storage but no extra applications of the operator or preconditioner,
and, we minimize
the unpreconditioned residual $2$-norm rather than in some other norm; 
see \cite{SS.2007} for
more details.

\section{Preliminaries}\label{section.prelim}
We begin with a brief review of Krylov subspace methods as well as
techniques of subspace recycling and for solving
shifted linear system.
Recall that in many Krylov subspace iterative methods for solving
\be\label{eqn.Axb}
	\vek A\vek x = \vek b
\ee
with $\vek A\in\Cnn$ , we generate an orthonormal basis for 
\be\label{eqn.Krylov-basis}
\CK_{j}(\vek A,\vek u) = \text{span}\curl{\vek u, \vek A\vek u,\ldots, \vek A^{j-1}\vek u}
\ee
with the Arnoldi process, where $\vek u$ is some starting vector.  Let $\vek V_{j}\in\C^{n\times j}$ be the matrix with orthonormal columns generated by the Arnoldi process spanning $\CK_{j}(\vek A,\vek u)$.  Then we have the Arnoldi relation
\be\label{eqn.arnoldi-relation}
	\vek A\vek V_{j} = \vek V_{j+1}\overline{\vek H}_{j}
\ee
with $\overline{\vek H}_{j}\in\C^{(j+1)\times j}$; see, e.g., \cite[Section 6.3]{Saad.Iter.Meth.Sparse.2003} and \cite{szyld.simoncini.survey.2007}.  Let $\vek x_{0}$ be an initial approximation to the solution of a linear system we wish to solve and $\vek r_{0}=\vek b-\vek A\vek x_{0}$ be the initial residual. At iteration $j$, we choose 
$\vek x_{j}=\vek x_{0}+\vek t_{j}$, with $\vek t_{j}\in\CK_{j}(\vek A,\vek r_{0})$.  In 
GMRES \cite{Saad.GMRES.1986}, $\vek t_{j}$ satisfies
\[
	\vek t_{j} = \argmin{\vek t\in \CK_{j}(\vek A,\vek r_{0})}\norm{\vek b-\vek A(\vek x_{0} + \vek t)},
\]
which is equivalent to solving the smaller minimization problem
\be\label{eqn.gmres-min}
	\vek y_{j} = \argmin{\vek y\in\C^{j}}\norm{ \overline{\vek H}_{j}\vek y - \norm{\vek r_{0}}\vek e_{1}^{(j+1)}},
\ee
where $\vek e_{J}^{(i)}$ denotes the $J$th Cartesian basis vector in $\C^{i}$.
We then set $\vek x_{j} = \vek x_{0} + \vek V_{j}\vek y_{j}$. Recall that in restarted 
GMRES, often called GMRES($m$), we run an $m$-step cycle of the GMRES method and 
compute an approximation $\vek x_{m}$.   We halt the process, discard $\vek V_{m}$, 
and restart with the new residual.  This process is repeated until we 
achieve convergence.  An adaption of restarted GMRES to solve (\ref{eqn.one-shifted-fam}) has 
been previously proposed; see, e.g., \cite{Frommer1998}.  

Many methods for the simultaneous solution of shifted systems 
(see, e.g., \cite{Darnell2008,Frommer2003,Frommer1998,Frommer1995,Kirchner2011,Simoncini2003a}) 
take advantage of the fact that for any shift $\sigma\in\C$, the Krylov subspace generated by $\vek A$ and $\vek b$ is invariant under the shift, i.e., 
\be\label{eqn.shift-invariance}
	\CK_{j}(\vek A,\vek b) = \CK_{j}(\vek A+\sigma\vek I,\tilde{\vek b}),
\ee 
as long as the starting vectors are collinear, i.e., $\tilde{\vek b} = \beta\vek b$ for some 
$\beta\in\C\setminus\curl{0}$, with a shifted Arnoldi 
relation similar to~(\ref{eqn.arnoldi-relation})
\be\label{eqn.shifted-arnoldi-relation}
	(\vek A+\sigma\vek I)\vek V_{j} = \vek V_{j+1}\overline{\vek H}_{j}\Sfs,
\ee
where 
$\overline{\vek H}_{j}\Sfs = \overline{\vek H}_{j} + \sigma\begin{bmatrix}\vek I_{m
\times m}\\\vek 0_{1\times m}\end{bmatrix}$.
This collinearity must be maintained at restart. In \cite{SSX.2014}, this was shown to be a
troublesome restriction when attempting to extend such techniques augmented Krylov methods.
In the case of GMRES, Frommer and Gl\"assner were able to overcome this by minimizing only one residual
in the common Krylov subspace and forcing the others to be collinear.  This strategy also works in the case of
GMRES with deflated restarts \cite{Darnell2008} because of properties of the augmented space generated using harmonic
Ritz vectors.  However, it was shown
in \cite{SSX.2014} that residual collinearity cannot be enforced in general.
Furthermore, it is not compatible with general 
preconditioning. The invariance
(\ref{eqn.shift-invariance}) 
can lead to great savings
in memory costs; but with a loss of algorithmic flexibility.  Thus in Section \ref{section.dir-proj-framework}, 
we explore an alternative.

We briefly review \RGMRES\ for non-Hermitian $\vek A$.
Augmentation techniques designed specifically for Hermitian
linear systems have also been proposed; see, e.g., 
\cite{KMR.2001,SYEG.2000,WSG.2007}.  For a more general framework
for these types of methods, see \cite{GGL.2013}, elements of which form a part of the
thesis of Gaul \cite{Gaul.2014-phd}, which contains a wealth of information on this
topic.  Gaul and
Schl\"{o}mer describe recycling techniques in the context of 
self-adjoint operator equations in a general Hilbert space 
\cite{GS.2013-arXiv}.

We begin by clarifying what we mean by Recycled GMRES. We use this expression to describe the general category
of augmented GMRES-type methods which are then differentiated by the choice of augmenting subspace.
As we subsequently explain, these methods can all be formulated as a GMRES iteration being applied to a 
linear system premultiplied with a projector.  The intermediate solution to this projected 
problem can then be further
corrected yielding a minimum residual approximation for the original problem over an augmented Krylov subspace.  
GCRO-DR \cite{Parks.deSturler.GCRODR.2005} is one such method in this category, in which the 
augmented subspace is built from harmonic Ritz vectors.

The GCRO-DR method represents the confluence of two approaches: those descending from 
the implicitly restarted Arnoldi method \cite{LS.1996}, 
such as Morgan's GMRES-DR \cite{Morgan.GMRESDR.2002}, 
and those descending from de Sturler's GCRO method \cite{deSturler.GCRO.1996}. 
GMRES-DR is a restarted GMRES algorithm, where at the end of each cycle, harmonic Ritz vectors
are computed, and a subset of them is used to augment the Krylov subspace generated at the next cycle.
The GCRO method allows the user to select the optimal correction over arbitrary subspaces.  
This concept is extended by de Sturler in \cite{deSturler.GCROT.1999}, where 
a framework is provided to optimally reduce convergence rate slowdown due to discarding information upon restart.  
This algorithm is called GCROT, where OT stands for optimal truncation.  A simplified version of the GCROT approach, 
based on restarted GMRES (called LGMRES) is presented in \cite{Manteuffel.LGMRES.2005}. 
Parks et~al.\ in
\cite{Parks.deSturler.GCRODR.2005} combine the ideas of  \cite{Morgan.GMRESDR.2002} 
and \cite{deSturler.GCROT.1999} 
and extend them to a sequence of slowly-changing linear systems. They call their method GCRO-DR.
This method and GCROT are Recycled GMRES methods. 

Suppose we are solving (\ref{eqn.Axb}), and we have a $k$-dimensional subspace $\CU$
whose image under the action of $\vek A$ is $\CC = \vek A\,\CU$.  Let $\vek{P}$ be the 
orthogonal projector onto $\CC^{\perp}$.  Let $\vek x_{0}$ be such that
$\vek r_{0}\in\CC^{\perp}$.  At iteration $m$, 
the \RGMRES\ method generates the approximation 
\be\nn
	\vek x_{m}=\vek x_{0}+\vek s_{m} + \vek t_{m}
\ee
where $\vek s_{m}\in\CU$ and $\vek t_{m}\in\CK_{m}(\vek{P}\vek A,\vek r_{0})$.  The corrections
$\vek s_{m}$ and $\vek t_{m}$ are chosen according to the minimum residual, Petrov-Galerkin
condition over the augmented Krylov subspace, i.e.,
\be\label{eqn.rgmres-petrov-galerkin}
	\vek r_{m}\perp \vek A\prn{\CU + \CK_{m}\prn{\vek{P}\vek A,\vek r_{0}}}.  
\ee
At the end of the cycle,
an updated $\CU$ is constructed, the Krylov subspace basis is discarded, and we restart.
At convergence, $\CU$ is saved, to be used when solving the next linear system.
This process is equivalent to applying GMRES to the projected problem
\be\label{eqn.rGMRES-proj-prob}
	\vek P\vek A\prn{\widehat{\vek x}_{0} + \vek t} = \vek P\vek b
\ee
where $\vek t_{m}$ is the $m$th GMRES correction for \eqref{eqn.rGMRES-proj-prob}
the second correction $\vek s_{m}\in\CU$ is the orthogonal projection of 
$\vek t_{m}$ onto $\CU$ where the orthogonality is with respect to the inner product
induced by the positive-definite matrix $\vek A^{\ast}\vek A$
\footnote[2]{We can write explicitly $\vek s_{m} = \vek P_{\CU}\vek t_{m}$ where we define
$\vek P_{\CU} = \vek U\prn{\vek U^{\ast}\vek A^{\ast}\vek A\vek U}^{-1}\vek U^{\ast}\vek A^{\ast}\vek A$ which can be rewritten $\vek P_{\CU}=\vek U\vek C^{\ast}\vek A$}; see, e.g., \cite{Gaul.2014-phd,GGL.2013}.

\RGMRES\ can be described as a modified GMRES iteration.  
Let $\vek U\in\C^{n\times k}$ have columns spanning $\CU$, scaled such that
$ 
	\vek C=\vek A\vek U
$ 
has orthonormal columns.  Then we can apply
$\vek{P} = \vek I - \vek C\vek C^{\ast}$ to $\vek A\vek v_{j}$ using 
$k$ steps of the Modified Gram-Schmidt process. 
The orthogonalization coefficients
are stored in the $m$th column of $\vek B_{m} = \vek C^{\ast}\vek A\vek V_{m}$,
which is simply $\vek B_{m-1}$ with one new column appended.
Let $\overline{\vek H}_{m}$ and $\vek V_{m}$ be defined as before, but for the projected
Krylov subspace $\CK_{m}\prn{\vek{P}\vek A,\vek r_{0}}$.  Enforcing (\ref{eqn.rgmres-petrov-galerkin})
is equivalent to solving the GMRES minimization problem (\ref{eqn.gmres-min}) for
 $\CK_{m}\prn{\vek{P}\vek A,\vek r_{0}}$ and setting 
 \be\nn
	\vek s_{m}=-\vek U\vek B_{m}\vek y_{m}\mbox{\ \ and\ \ }\vek t_{m}=\vek V_{m}\vek y_{m},
\ee
so that 
\begin{align} \nn
	\vek x_{m}& = \vek x_{0} - \vek U\vek B_{m}\vek y_{m} + \vek V_{m}\vek y_{m}
= \vek x_{0} + \begin{bmatrix} \vek U & \vek V_{m}  \end{bmatrix}
			       \begin{bmatrix}  -\vek B_{m}\vek y_{m}\\ \vek y_{m}\end{bmatrix} .
\end{align}
This is a consequence of the fact that the \RGMRES\ least squares problem, as stated in \cite[Equation 2.13]{Parks.deSturler.GCRODR.2005} can be satisfied exactly in the first $k$ rows, and this was first observed in 
\cite{deSturler.GCRO.1996}.  The choice of the subspace $\CU$ then determines the actual method.

\section{A direct projection framework}\label{section.dir-proj-framework}
We develop a general framework of minimum residual methods for shifted
linear systems which encompasses both unpreconditioned and preconditioned systems.
We propose to solve both a single family of shifted systems \eqref{eqn.one-shifted-fam} and
sequences of shifted system families of the form \eqref{eqn.one-shifted-family-seq}.
However, it suffices to propose our method in a simpler setting in which we drop the index $i$ and assume there are only
two systems, a base system and a shifted system. 
Thus for simplicity, we restrict our description to two model problems:
the unpreconditioned problem
\be\label{eqn.model-problem-unprec}
	\vek A\vek x = \vek b \mbox{\ \ and\ \ }(\vek A + \sigma\vek I)\vek x\Sfs = \vek b
\ee
and the right-preconditioned problem
\be\label{eqn.model-problem-prec}
	\vek A\vek M^{-1}\vek w = \vek b \mbox{\ \ and\ \ }(\vek A + \sigma\vek I)\vek M^{-1}\vek w\Sfs = \vek b
\ee
where $\vek w_{0} = \vek M\vek x_{0}$ and $\vek w_{0}\Sfs = \vek M\vek x_{0}\Sfs$, and
after $m$ iterations we set $\vek x_{m} = \vek M^{-1}\vek w_{m}$ and 
we set $\vek x_{m}\Sfs = \vek M^{-1}\vek w_{m}\Sfs$.  
In this setting, we can propose minimum residual Krylov subspace methods in the cases that we do and do not
have an augmenting subspace $\CU$.

We describe the proposed methods in terms of a general sequence of nested subspaces
\be\nn
	\CS_{1}\subset\CS_{2}\subset\cdots\CS_{m}\subset\cdots
\ee
This allows us to cleanly present these techniques as minimum residual 
projection methods and later to give clear analysis, applicable to any method
fitting into this framework.  Then we can derive different methods by specifying
$\CS_{m}$, e.g., \linebreak $\CS_{m}=\CK_{m}(\vek A,\vek r_{0})$.

Let $\curl{\CS_{m}}_{i=1}^{m}$ be the nested sequence of subspaces produced by some
some iterative method for solving \eqref{eqn.model-problem-unprec} or 
\eqref{eqn.model-problem-prec}, after $m$ iterations.  
In the unpreconditioned case \eqref{eqn.model-problem-unprec},
suppose we have initial approximations $\vek x_{0}$ and $\vek x_{0}\Sfs$ for the base and shifted 
systems, respectively.  For conciseness, let us denote $\vek A\Sfs = \vek A + \sigma \vek I$.
At iteration $m$, we compute corrections $\vek t_{m},\,\vek t_{m}\Sfs\in\CS_{m}$ 
which satisfy the minimum residual conditions
\be\label{eqn.unprec-Galerkin}
	\vek b - \vek A(\vek x_{0} + \vek t_{m})\perp \vek A\CS_{m}\mbox{\ \ and\ \ }\vek b - \vek A\Sfs(\vek x_{0}\Sfs + \vek t_{m}\Sfs)\perp \vek A\Sfs\CS_{m}.
\ee
In the preconditioned case \eqref{eqn.model-problem-prec},
suppose we begin with initial approximations 
$\vek w_{0}=\vek M\vek x_{0}$ and 
$\vek w_{0}\Sfs=\vek M\vek x_{0}\Sfs$.
Let us denote the preconditioned operators
\be\nn
\vek A_{p} = \vek A\vek M^{-1}\mbox{\ \ and\ \ }\vek A_{p}\Sfs = (\vek A + \sigma\vek I)\vek M^{-1}
\ee
At iteration $m$, we compute corrections $\vek t_{m},\,\vek t_{m}\Sfs\in\CS_{m}$ 
which satisfy the minimum residual conditions
\be\label{eqn.prec-Galerkin}
	\vek b - \vek A_{p}(\vek w_{0} + \vek t_{m})\perp \vek A_{p}\CS_{m}\mbox{\ \ and\ \ }\vek b - \vek A_{p}\Sfs (\vek w_{0}\Sfs + \vek t_{m}\Sfs)\perp \vek A_{p}\Sfs \CS_{m}.
\ee
We emphasize that \textit{the same preconditioner} is used for all systems.  

In this framework, we assume that the minimizer for the base case is constructed via a
predefined iterative method, the method which generates the sequence $\curl{\CS_{m}}$.  Therefore,
it suffices to describe the residual projection for the shifted system.  We can write the update
of the shifted system approximation by explicitly constructing the orthogonal projector which
is applied during a Petrov-Galerkin projection.  Let $\curl{\vek s_{1},\vek s_{2},\ldots,\vek s_{m}}$
be a basis for $\CS_{m}$ which we take as the columns of $\vek S_{m}\in\C^{n\times m}$.  
Then we can write this projection and update
\begin{align}
	\vek r_{m}\Sfs &= \vek r_{0}\Sfs - {\vek A\Sfs\vek S_{m}}\vek y_{m}\Sfs\mbox{\ \ and}\nn\\
	\vek x_{m}\Sfs &= \vek x_{0}\Sfs + {\vek S_{m}}\vek y_{m}\Sfs\label{eqn.GMRES-shift-proj-unprec}
\end{align}
where $\vek w_{m}\Sfs = \vek N_{m}\Sfs^{-1}\prn{\vek A\Sfs\vek S_{m}}^{\ast}\vek r_{0}\Sfs$ and 
$\vek N_{m}\Sfs = \prn{\vek S_{m}^{\ast}\vek A\Sfs^{\ast}\vek A\Sfs\vek S_{m}}$ is the projection scaling matrix, 
since we assume that 
$\vek A\Sfs\vek S_{m}$ does not have orthonormal columns.  
For well-chosen $\CS_{m}$, these projections can be applied using 
already-computed quantities.

In the following subsections, we derive new methods by specifying
subspaces $\curl{\CS_{m}}$ and a matrix $\vek S_{m}$.  
This will define $\vek N_{m}\Sfs$.  
We show that for these choices, $\vek N_{m}\Sfs$ is composed of blocks
which can be built from already-computed quantities.  Thus, for appropriate
choices of $\CS_{m}$, either \eqref{eqn.unprec-Galerkin}
or \eqref{eqn.prec-Galerkin}, can be applied with manageable additional costs. 

We highlight that a strength of this framework that we can 
develop methods for shifted systems on top of an existing
iterative methods, with a few modifications.  As the framework only requires
a sequence of nested subspaces, it is completely compatible with
with both standard Krylov subspace methods as well as  
flexible and inexact Krylov subspace methods.

\subsection{A GMRES method for shifted systems}\label{subsection.GMRES-shift}
In the case that we apply the GMRES iteration to the base system, at iteration $m$, our search
space is $\CS_{m}:=\CK_{m}(\vek A, \vek r_{0})$, and
the matrix $\vek S_{m}:=\vek V_{m}$ has the first $m$ Arnoldi 
vectors as columns.  The projection and update  
\eqref{eqn.GMRES-shift-proj-unprec} can be simplified due to the shifted Arnoldi
relation \eqref{eqn.shifted-arnoldi-relation}.  The matrix 
$\vek N_{m}\Sfs:=\overline{\vek H}_{m}\Sfs^{\ast}\overline{\vek H}_{m}\Sfs\in\C^{m\times m}$ can be constructed
from the already computed upper Hessenberg matrix. Thus the projection 
\eqref{eqn.unprec-Galerkin} can be rewritten 
\begin{align}
	\vek x_{m}\Sfs &= \vek x_{0}\Sfs + {\vek V_{m}\Sfs}\vek y_{m}\Sfs\mand\nn\\
	\vek r_{m}\Sfs &= \vek r_{0}\Sfs - {\vek V_{m+1}\overline{\vek H}_{m}\Sfs}\vek y_{m}\Sfs\nn
\end{align}
where $\vek y_{m}\Sfs = \prn{\overline{\vek H}\Sfs^{\ast}\overline{\vek H}\Sfs}^{-1}{{\overline{\vek H}_{m}\Sfs}}^{\ast}\vek V_{m+1}^{\ast}\vek r_{0}\Sfs$.
As it can be appreciated, applying this is equivalent to solving the least squares
problem
\be\label{eqn.y-shift-LS}
	\vek y_{m}\Sfs = \argmin{\vek y\in\C^{i}}\norm{\overline{\vek H}_{m}\Sfs\vek y - \vek V_{m+1}^{\ast}\vek r_{0}\Sfs}
\ee
and setting $\vek x_{m}\Sfs = \vek x_{0}\Sfs + \vek V_{m}\vek y_{m}\Sfs$.
This method has similarities with the GMRES method for shifted systems of 
Frommer and Gl\"{a}ssner \cite{Frommer1998}, which is derived from the invariance
\eqref{eqn.shift-invariance}.  In the method proposed in \cite{Frommer1998}, one must solve small linear
systems for each shifted system whereas here one must solve the small least-squares problem \eqref{eqn.y-shift-LS}.
The main difference is that what we propose does not guarantee convergence of all system in one Krylov subspace
whereas in \cite{Frommer1998}, this is guaranteed under certain conditions.  The strength here comes from the ability
to precondition.
\subsubsection{Preconditioning}
Introducing preconditioning into this setting presents complications.  No longer can we use 
the shifted Arnoldi relation \eqref{eqn.shifted-arnoldi-relation} as we could in the 
unpreconditioned case.  However, by storing some extra vectors, as in Flexible GMRES \cite{Saad1993}, 
one can enforce \eqref{eqn.prec-Galerkin} with no additional application of the operator or 
preconditioner.  

Recall that in right-preconditioned GMRES (see, e.g., 
\cite[Sections 9.3.2 and 9.4.1]{Saad.Iter.Meth.Sparse.2003}) that 
$\CS_{m}:=\vek M^{-1}\CK(\vek A_{p},\vek r_{0})$,
and $
\vek S_{m}:=\vek M^{-1}\vek V_{m}$.  This space
is never explicitly constructed, though, since if $\vek y_{m}$ is the solution to the 
GMRES least squares problem \eqref{eqn.gmres-min} in the preconditioned
case, we simply set $\vek x_{m} = \vek x_{0} + \vek M^{-1}\prn{\vek V_{m}\vek y_{m}}$.
However, in flexible GMRES, one must store this basis.  
For all $1\leq i\leq m$, let 
$\vek z_{i} = \vek M^{-1}\vek v_{i}$, and let these vectors be the columns of 
$\vek Z_{m}\in\C^{n\times i}$ so that $\vek Z_{m} = \vek M^{-1}\vek V_{m}$.  

With these 
vectors, one can enforce \eqref{eqn.prec-Galerkin}. 
Observe that we can write $\vek A_{p}\Sfs = \vek A_{p} + \sigma\vek M^{-1}$.
We explicitly project the residual, but this time
 onto $\curl{\vek A_{p}\Sfs\CK_{m}(\vek A_{p},\vek r_{0})}^{\perp}$,
 \be\label{eqn.shift-gmres-prec-proj}
 	\vek r_{m}\Sfs = \vek r_{0}\Sfs - \prn{\vek A_{p} + \sigma\vek M^{-1}}\vek V_{m}\vek N_{m}\Sfs^{-1}\brac{\prn{\vek A_{p} + \sigma\vek M^{-1}}\vek V_{m}}^{\ast}\vek r_{0}\Sfs
 \ee
 where $\vek N_{m}\Sfs =\brac{\prn{\vek A_{p} + \sigma\vek M^{-1}}\vek V_{m}}^{\ast}\brac{\prn{\vek A_{p} + \sigma\vek M^{-1}}\vek V_{m}}$.  With the right-preconditioned shifted Arnoldi relation
 \be\nn
 	\prn{\vek A_{p} + \sigma\vek M^{-1}}\vek V_{m} = \vek V_{m+1}\overline{\vek H}_{m} + \sigma \vek Z_{m}
 \ee
we rewrite 
 \be\nn
 	\vek N_{m}\Sfs = \overline{\vek H}_{m}^{\ast}\overline{\vek H}_{m}+ \sigma\overline{\vek H}_{m}^{\ast}\vek V_{m+1}^{\ast}\vek Z_{m} + \overline{\sigma}\vek Z_{m}^{\ast}\vek V_{m+1}\overline{\vek H}_{m} + \ab{\sigma}^{2}\vek Z_{m}\vek Z_{m}.
 \ee
 Thus, the approximation update and
 the residual projection \eqref{eqn.shift-gmres-prec-proj} can be rewritten
 \begin{align}
 	\vek x_{m}\Sfs &= \vek x_{0}\Sfs + \vek Z_{m}\vek y_{m}\Sfs\nn\\
 	\vek r_{m}\Sfs &= \vek r_{0}\Sfs - \prn{\vek V_{m+1}\overline{\vek H}_{m} + \sigma \vek Z_{m}}\vek y_{m}\Sfs.\nn
 \end{align}
 where $\vek y_{m}\Sfs = \vek N_{m}\Sfs^{-1}\brac{\prn{\vek V_{m+1}\overline{\vek H}_{m} + \sigma \vek Z_{m}}}^{\ast}\vek r_{0}\Sfs$.
 This projection process involves only the precomputed matrices ($\overline{\vek H}_{m}$, 
 $\vek V_{m+1}$, and $\vek Z_{m+1}$).
 The matrices $\overline{\vek H}_{m}^{\ast}\overline{\vek H}_{m} $, 
 $\overline{\vek H}_{m}^{\ast}\vek V_{m+1}^{\ast}\vek Z_{m}$, and $\vek Z_{m}^{\ast}\vek Z_{m}$
 can be computed once, independent of the number of shifted systems.  The solution of a 
 dense Hermitian linear system with 
 $\vek N_{m}\Sfs$ must be performed for each $\sigma$.  This solution
 of a Hermitian $m\times m$ linear system costs 
 $\CO(m^{3})$ floating point operations (FLOPS).  The right-preconditioned shifted GMRES
 algorithm (sGMRES)
 is shown in Algorithm \ref{algorithm.sgmres}.  Observe
 that an implementation can rely heavily on an existing GMRES code.
 It should be noted that all but one step of the 
 shifted residual projections can be formulated
 in terms of block/BLAS-3 operations so that most computations for all
 shifts are performed simultaneously.  
 \input{sgmres-algorithm.tex} 
\subsection{An rGMRES method for shifted systems}\label{subsection.rgmres-shift}
Suppose now that our iteration for the base system is a \RGMRES\ method.  

We begin by projecting the initial residual 
$\vek r_{-1}\Sfs$
associated to initial approximation $\vek x_{-1}\Sfs$, so that we begin with 
$\vek r_{0}\Sfs\perp \vek A\Sfs\,\CU$.  
This is equivalent to computing the minimum residual correction $\vek t_{0}\Sfs\in\CU$ and setting $\vek x_{0}\Sfs = \vek x_{-1}\Sfs + \vek t_{0}\Sfs$.
In \RGMRES, such a projection is
necessary to correctly derive the algorithm.  
For the shifted system, the projection is not necessary,
but it does allow for some simplifications later in the derivation.  We have then, 
\be\label{eqn.init-recyc-proj}
	\vek x_{0}\Sfs = \vek x_{-1}\Sfs + \vek U\vek y_{0}\Sfs\mbox{\ \ and\ \ }\vek r_{0}\Sfs = \vek r_{-1}\Sfs - {\vek A\Sfs\vek U}\vek y_{0}\Sfs,
\ee
where $\vek y_{0}\Sfs = \vek N_{0}\Sfs^{-1}\prn{\vek A\Sfs\vek U}^{\ast}\vek r_{-1}$ and $\vek N_{0}\Sfs = \prn{\vek A\Sfs\vek U}^{\ast}\prn{\vek A\Sfs\vek U}$.  
\linebreak Since $\vek A\Sfs\vek U = \vek C + \sigma\vek U$, this projection can 
be simplified and computed with manageable additional expense,
\be\nn
	\vek r_{0} = \vek r_{-1} - (\vek C+\sigma\vek U)\vek N_{m}\Sfs^{-1}(\vek C+\sigma\vek U)^{\ast}\vek r_{-1}
\ee
where we rewrite $\vek N_{0}\Sfs = \vek I_{k\times k} + \sigma\vek C^{\ast}\vek U + \overline{\sigma}\vek U^{\ast}\vek C + \ab{\sigma}^{2}\vek U^{\ast}\vek U$.  The matrices $\vek C^{\ast}\vek U$ and 
$\vek U^{\ast}\vek U$ must only be computed once, regardless of the number of shifts, and for each shift we solve 
$\vek N_{0}\Sfs\vek y_{0}\Sfs = (\vek C + \sigma\vek U)^{\ast}\vek r_{-1}\Sfs$.

After a cycle of \RGMRES\ for the base system, \eqref{eqn.unprec-Galerkin}
must be enforced for each shifted system.  
At iteration $m$, our search space 
\linebreak $\CS_{m}\:= \CU + \CK_{m}(\vek P\vek A_{p}, \vek r_{0})$.
The augmented matrix $\vek S_{m}:= \begin{bmatrix} \vek U & \vek V_{m}\end{bmatrix}$
contains as columns the basis for $\CU$ and $\CK_{m}(\vek P\vek A_{p}, \vek r_{0})$.
In this case, we have $\vek N_{m}\Sfs = \curl{(\vek A+\sigma\vek I)\begin{bmatrix}\vek U & \vek V_{m}\end{bmatrix} }^{\ast}\curl{(\vek A+\sigma\vek I)\begin{bmatrix}\vek U & \vek V_{m}\end{bmatrix} }$.
From \cite{SSX.2014}, we have the identity
\be\nn
	(\vek A+\sigma\vek I)\begin{bmatrix}\vek U & \vek V_{m}\end{bmatrix} = \begin{bmatrix} \prn{\vek C+\sigma\vek U} & \prn{\vek C\vek B_{m}+\vek V_{m+1}\overline{\vek H}_{m} + \sigma\vek V_{m}} \end{bmatrix}.
\ee
Thus, in the unpreconditioned case, for the augmented Krylov subspace, we
can rewrite \eqref{eqn.unprec-Galerkin}
\begin{align}
	\vek r_{m}\Sfs &= \vek r_{0}\Sfs - \begin{bmatrix} \prn{\vek C+\sigma\vek U} & \prn{\vek C\vek B_{m}+\vek V_{m+1}\overline{\vek H}_{m} + \sigma\vek V_{m}} \end{bmatrix}\vek y_{m}\Sfs\mbox{\ \ and}\nn\\
	\vek x_{m}\Sfs &= \vek x_{0}\Sfs + \begin{bmatrix}\vek U & \vek V_{m}\end{bmatrix}\vek y_{m}\Sfs
\end{align}
where $\vek y_{m}\Sfs = \vek N_{m}\Sfs^{-1}\begin{bmatrix} \prn{\vek C+\sigma\vek U} & \prn{\vek C\vek B_{m}+\vek V_{m+1}\overline{\vek H}_{m}^{(\sigma)}} \end{bmatrix}^{\ast}\vek r_{0}\Sfs$ and
\be\nn
	\resizebox{.9\hsize}{!}{$\vek N_{m}\Sfs = \begin{bmatrix}
										\vek I + \sigma\vek C^{\ast}\vek U+ \overline{\sigma}\vek U^{\ast}\vek C  + \ab{\sigma}^{2}\vek U^{\ast}\vek U & \vek B_{m} + \overline{\sigma}\vek U^{\ast}\vek C\vek B_{m} + \overline{\sigma}\vek U^{\ast}\vek V_{m+1}\overline{\vek H}_{m} + \ab{\sigma}^{2}\vek U^{\ast}\vek V_{m}\\
										\vek B_{m}^{\ast} + \sigma\vek B_{m}^{\ast}\vek C^{\ast}\vek U + \sigma \overline{\vek H}_{m}^{\ast}\vek V_{m+1}^{\ast}\vek U  + \ab{\sigma}^{2}\vek V_{m}^{\ast}\vek U& \vek B_{m}^{\ast}\vek B_{m} + \overline{\vek H}_{m}^{\ast}\overline{\vek H}_{m} + {\sigma\vek H_{m} + \overline{\sigma}\vek H_{m}^{\ast}} + \ab{\sigma}^{2}\vek I
								   \end{bmatrix}.$}
\ee
This projection can be performed using 
already computed quantities, and the matrices
$\vek U^{\ast}\vek C$, $\vek U^{\ast}\vek U$, $\vek U^{\ast}\vek C\vek B_{m}$, 
$\vek U^{\ast}\vek V_{m+1}\overline{\vek H}_{m}$, $\overline{\vek H}_{m}^{\ast}\overline{\vek H}_{m}$, $\vek H_{m}$, and 
$\vek B_{m}^{\ast}\vek B_{m} $ need only be computed once, regardless of the 
number of shifts.  The computation
of $\vek y_{m}\Sfs$ must be performed for every shift at a cost of $\CO((m+k)^{3})$.

\subsubsection{Preconditioning}

Introducing right preconditioning creates some difficulties which we can
again surmount by storing some extra vectors.  
We note that in the case of preconditioning, we have $\vek C = \vek A\vek M^{-1}\vek U$.
In this case, for right
preconditioned \RGMRES, the search space for the base system is 
$\CS_{m}:=\vek M^{-1}\curl{\CU + \CK_{m}(\vek P\vek A,\vek r_{0})}$.
Let 
$\vek Z_{\CU} = \vek M^{-1}\vek U$ and $\vek Z_{m}=\vek M^{-1}\vek V_{m}$, as in Section \ref{subsection.GMRES-shift}. 

Using $\vek Z_{\CU}$, we can cheaply perform the initial residual projection,
\begin{align}
	\vek x_{0}\Sfs &= \vek x_{-1}\Sfs + \vek U\vek y_{0}\Sfs\mbox{\ \ and}\nn\\
	\vek r_{0}\Sfs & =\vek r_{-1}\Sfs - \prn{\vek A_{p}\Sfs\vek U}\vek y_{0}\Sfs\label{eqn.init-recyc-proj-prec}
\end{align}
where $\vek y_{0}\Sfs=\vek N_{0}\Sfs^{-1}\prn{\vek A_{p}\Sfs\vek U}^{\ast}\vek r_{-1}\Sfs$ and
$\vek N_{0}\Sfs = \prn{\vek A_{p}\Sfs\vek U}^{\ast}\prn{\vek A_{p}\Sfs\vek U}$.
We can write 
\be\nn
	\vek A_{p}\Sfs\vek U = \vek C + \sigma\vek Z_{\CU}.
\ee
The subspace $\CU$ either is available from at the start of the 
algorithm (in which case $\vek U$ must be scaled so that 
$\vek A_{p}\vek U=\vek C$ has orthonormal columns), or it is 
constructed at the end of a
restart cycle.  In either case, $\vek Z_{\CU}$ is available in the course
of the computation and can be saved.  Thus the projection
\eqref{eqn.init-recyc-proj} can be performed with already computed
quantities,
\begin{align}
	\vek x_{0}\Sfs &= \vek x_{-1}\Sfs + \vek U\vek y_{0}\Sfs\mbox{\ \ and}\nn\\ \vek r_{0}\Sfs &= \vek r_{-1}\Sfs - \prn{\vek C + \sigma\vek Z_{\CU}}\vek y_{0}\Sfs\label{eqn.init-recyc-proj-prec-simpl},
\end{align}
where we rewrite $\vek y_{0}\Sfs = \vek N_{0}\Sfs^{-1}\prn{\vek C + \sigma\vek Z_{\CU}}^{\ast}\vek r_{-1}\SfAlg{\sigma_{\ell}}$ and

$$\vek N_{0}\Sfs = \vek I + \sigma\vek C^{\ast}\vek Z_{\CU} + \overline{\sigma}\vek Z_{\CU}^{\ast}\vek C + \ab{\sigma}^{2}\vek Z_{\CU}^{\ast}\vek Z_{\CU}.$$

After a cycle of right-preconditioned \RGMRES, we must perform the projection
\eqref{eqn.prec-Galerkin} for each shifted system.  
We proceed slightly differently in this derivation than in 
the unpreconditioned case.  We have 
\be\nn
\vek N_{m}\Sfs:=	\curl{\vek A_{p}\Sfs \begin{bmatrix}\vek U & \vek V_{m}\end{bmatrix}}^{\ast}\curl{\vek A_{p}\Sfs \begin{bmatrix}\vek U & \vek V_{m}\end{bmatrix}} 
\ee
Following from  \cite{Parks.deSturler.GCRODR.2005}, we define 
\be\nn
	\overline{\vek G}_{m} = \begin{bmatrix}
												\vek I_{k\times k} & \vek B_{m}\\
												\vek 0_{(m+1)\times k} & \overline{\vek H}_{m}	
										  \end{bmatrix},
\ee
which yields the augmented Arnoldi relation
\be\label{eqn.aug-arn-rel}
	\vek A_{p}\begin{bmatrix}\vek U & \vek V_{m}\end{bmatrix} = \begin{bmatrix}\vek C & \vek V_{m+1}\end{bmatrix}\overline{\vek G}_{m}.
\ee
Using the relation \eqref{eqn.aug-arn-rel}, an identity 
for the shifted operator with right preconditioning follows,
\be\label{eqn.shifted-prec-aug-arnoldi}
	\vek A_{p}\Sfs\begin{bmatrix}\vek U & \vek V_{m}\end{bmatrix} = \begin{bmatrix}\vek C & \vek V_{m+1}\end{bmatrix}\overline{\vek G}_{m} + \sigma \begin{bmatrix}\vek Z_{\CU} & \vek Z_{m}\end{bmatrix}.
\ee
We use the relation \eqref{eqn.shifted-prec-aug-arnoldi} to derive the
expansion
\begin{align}
	\vek N_{m}\Sfs& = \overline{\vek G}_{m}^{\ast}\overline{\vek G}_{m} + \ab{\sigma}^{2}\begin{bmatrix}\vek Z_{\CU}^{\ast}\vek Z_{\CU} & \vek Z_{\CU}^{\ast}\vek Z_{m}\\ \vek Z_{m}^{\ast}\vek Z_{\CU} & \vek Z_{m}^{\ast}\vek Z_{m}\end{bmatrix}	 
	+ \sigma\overline{\vek G}_{m}^{\ast}\begin{bmatrix}\vek C^{\ast}\vek Z_{\CU} & \vek C^{\ast}\vek Z_{m} \\ \vek V_{m+1}^{\ast}\vek Z_{\CU} & \vek V_{m+1}^{\ast}\vek Z_{m}\end{bmatrix}\nn \\
		&\quad + \overline{\sigma}\begin{bmatrix}\vek Z_{\CU}^{\ast}\vek C & \vek Z_{\CU}^{\ast}\vek V_{m+1} \\ \vek Z_{m}^{\ast}\vek C & \vek Z_{m}^{\ast}\vek V_{m+1}\end{bmatrix}\overline{\vek G}_{m}.\label{eqn.Nsig-prec-rgmres-sum}
\end{align}
Thus, the projection 
can be performed for
each shift using already computed quantities.  
This yields the following updates of the approximation and residual
 \input{sgcrodr-algorithm.tex} 
\begin{align}
	\vek x_{m}\Sfs &= \vek x_{0 }\Sfs + \begin{bmatrix}\vek Z_{\CU} & \vek Z_{m}\end{bmatrix}\vek y_{m}\Sfs \\
	\vek r_{m}\Sfs &= \vek r_{0}\Sfs - \curl{\begin{bmatrix}\vek C & \vek V_{m+1}\end{bmatrix}\overline{\vek G}_{m} + \sigma \begin{bmatrix}\vek Z_{\CU} & \vek Z_{m}\end{bmatrix}}\vek y_{m}\Sfs
\end{align}
where $\vek y_{m}\Sfs = \vek N_{m}\Sfs^{-1}\curl{\begin{bmatrix}\vek C & \vek V_{m+1}\end{bmatrix}\overline{\vek G}_{m} + \sigma \begin{bmatrix}\vek Z_{\CU} & \vek Z_{m}\end{bmatrix}}^{\ast}\vek r_{0}\Sfs$.
We observe that because of the initial projection of the shifted residual 
\eqref{eqn.init-recyc-proj-prec}, we can simplify
\begin{align}
	\resizebox{.45\hsize}{!}{$\curl{\begin{bmatrix}\vek C & \vek V_{m+1}\end{bmatrix}\overline{\vek G}_{m} + \sigma \begin{bmatrix}\vek Z_{\CU} & \vek Z_{m}\end{bmatrix}}^{\ast}\vek r_{0}\Sfs$} &=\resizebox{.45\hsize}{!}{$\curl{\begin{bmatrix} \vek C & \vek C\vek B_{m}+\vek V_{m+1}\overline{\vek H}_{m}\end{bmatrix} + \sigma \begin{bmatrix}\vek Z_{\CU} & \vek Z_{m}\end{bmatrix}}^{\ast}\vek r_{0}\Sfs$} \nn \\
	&= \begin{bmatrix}\vek C^{\ast}\vek r_{0}\Sfs\\ \vek B_{m}^{\ast}\vek C^{\ast}\vek r_{0}\Sfs+ \overline{\vek H}_{m}^{\ast}\vek V_{m+1}^{\ast}\vek r_{0}\Sfs\end{bmatrix} + \overline{\sigma}\begin{bmatrix}\vek Z_{\CU}^{\ast}\vek r_{0}\Sfs\\ \vek Z_{m}^{\ast}\vek r_{0}\Sfs\end{bmatrix}\nn \\
	&= \begin{bmatrix}\vek 0\\ \vek B_{m}^{\ast}\vek C^{\ast}\vek r_{0}\Sfs+ \overline{\vek H}_{m}^{\ast}\vek V_{m+1}^{\ast}\vek r_{0}\Sfs + \overline{\sigma}\vek Z_{m}^{\ast}\vek r_{0}\Sfs\end{bmatrix},\nn
\end{align}
and thus we can rewrite
\be\nn
	\vek y_{m}\Sfs = \vek N_{m}\Sfs^{-1}\begin{bmatrix}\vek 0\\ \vek B_{m}^{\ast}\vek C^{\ast}\vek r_{0}\Sfs+ \overline{\vek H}_{m}^{\ast}\vek V_{m+1}^{\ast}\vek r_{0}\Sfs + \overline{\sigma}\vek Z_{m}^{\ast}\vek r_{0}\Sfs\end{bmatrix}.
\ee
The matrices in the sum
\eqref{eqn.Nsig-prec-rgmres-sum} 
must be computed only once. For each shift, we must compute $\vek y_{m}\Sfs$ 
at a cost of $\CO\prn{(m+k)^{3}}$.   The right-preconditioned shifted \RGMRES\ 
 algorithm (srGMRES)
 is shown in Algorithm \ref{algorithm.sgcrodr}.  Observe
 that an implementation can rely heavily on an existing \RGMRES\ code.
 As in the case of Algorithm \ref{algorithm.sgmres}, all but one step of the 
 shifted residual projections can be formulated
 in terms of block/BLAS-3 operations so that almost all computations
 are performed simultaneously for all shifts.  
 We discuss costs further in Section \ref{section.cost}.
 
\section{Analysis of direct projection methods}\label{section.analysis}
In this section, we provide some analysis of the direction projection methods. 
We treat two issues in this section: quality of the approximations and cost of the methods.
\subsection{Quality of the approximations}\label{section.quality}

Since all residual corrections are minimum
residual projections, we can expect that, at worse, the projection of the
shifted residual will achieve no improvement. 

We follow the analysis presented in \cite{KMR.2001}.  This analysis follows from 
that presented in \cite{Chan1997} for case of Hermitian positive definite coefficient matrix.
In their analysis, 
the authors assume that a subset of eigenvectors (spanning subspace $\CY$) have been well-approximated
in the underlying Krylov subspace generated by QMR applied to
the base matrix (called the seed system in 
\cite{KMR.2001}).  The authors show that the performance of the QMR applied to 
the non-seed systems with projected residuals can 
be compared to that of a GMRES iteration in which $\CY$
has been projected away.  

In the case of Hermitian positive definite systems, analysis of the performance of CG-based seed-projection
was also extended to the case in which the coefficient matrix varies along with the right-hand side.  One of the
special cases considered is the present case, that one is solving a family of shifted linear 
systems.  In this case, one can again derive CG-based bounds dependent upon the set of eigenvectors
well-approximated by Ritz vectors generated by the CG iteration applied to the base system.

In extending this analysis, there are two complications. Algorithm \ref{algorithm.sgcrodr} does not minimize
over a Krylov subspace, and both methods may use preconditioning.  In either
case, we cannot easily leverage the polynomial approximation analysis.  
Also of concern is that GMRES-based methods applied to non-Hermitian problems of large dimension often
must be restarted, which does not need to be considered for the short-term recurrence-based methods
treated in \cite{CN.GalProjAnalMRHS.1999, Chan1997,KMR.2001}.  
However, if we restrict our analysis to Algorithm \ref{algorithm.sgmres}
without preconditioning (i.e., $\vek M=\vek I$) and do not consider restarting, we can analyze performance 
based on invariant subspace approximation. We follow from \cite{CN.GalProjAnalMRHS.1999,Chan1997} and specifically
use elements of analysis in \cite{KMR.2001} for the non-Hermitian case.

Let us assume that $\vek A$ is diagonalizable with eigendecomposition
\be\label{eqn.A-eig-decomp}
\vek A = \vek F\boldsymbol\Lambda\vek W\mwith\boldsymbol\Lambda = \diag\curl{\lambda_{1},\lambda_{2}\ldots,\lambda_{n}}\mand \vek W = \vek F^{-1},
\ee
with $\vek f_{i}$ being the $i$th column of $\vek F$ and $\vek w_{i}^{\ast}$ being the $i$th row
of $\vek W$. 
Consider the simplified problem \ref{eqn.model-problem-unprec},
where for the base system \eqref{eqn.Axb}, we have initial residual $\vek r_{0}$.
We first solve the base system using a GMRES iteration terminating in $j$ steps, generating the subspace 
$\CK_{j}(\vek A,\vek r_{0})$ with the associated $\vek V_{j}$ and $\overline{\vek H}_{j}$.  
Let $\vek x_{0}\Sfs$ be the initial approximation for the shifted system with residual 
$\vek r_{0}\Sfs = \vek b - \prn{\vek A+\sigma\vek I}\vek x_{0}\Sfs$.  Let 
$\hat{\vek r}_{0}\Sfs\perp\prn{\vek A+\sigma\vek I}\CK_{j}(\vek A,\vek r_{0})$ be the result of the
Lanczos-Galerkin projection of $\vek r_{0}\Sfs$ after the termination of GMRES applied to \eqref{eqn.Axb}.
If $\vek P$ is the projector onto $\CK_{j}(\vek A,\vek r_{0})$ which is orthogonal with respect to the
inner product induced by $\prn{\vek A+\vek \sigma\vek I}^{\ast}\prn{\vek A+\vek \sigma\vek I}$\footnote[2]{
i.e., $\vek P = \vek V_{m}\prn{\overline{\vek H}\Sfs^{\ast}\overline{\vek H}\Sfs}^{-1}{{\overline{\vek H}_{m}\Sfs}}^{\ast}\vek V_{m+1}^{\ast}$}, and
$\vek Q$ is the orthogonal projector onto $\prn{\vek A + \sigma\vek I}\CK_{j}(\vek A,\vek r_{0})$ with respect
to the Euclidean norm\footnote[3]{i.e., $\vek Q = {\vek V_{m+1}\overline{\vek H}_{m}\Sfs}\prn{\overline{\vek H}\Sfs^{\ast}\overline{\vek H}\Sfs}^{-1}{{\overline{\vek H}_{m}\Sfs}}^{\ast}\vek V_{m+1}^{\ast}$},
then we can then write the Lanczos-Galerkin projection as $\hat{\vek r}_{0}\Sfs = \prn{\vek I - \vek Q}\vek r_{0}\Sfs$,
and the associated updated approximation $\hat{\vek x}_{0}\Sfs$ results from the error projection 
\be\label{eqn.error-projection}
{\vek x\Sfs - \hat{\vek x}_{0}\Sfs} =  \prn{\vek I - \vek P}\prn{\vek x\Sfs - {\vek x}_{0}\Sfs}
\ee
This can be seen by studying the derivation in Section \ref{subsection.GMRES-shift} and is a general property
of minimum residual projections.  With this new starting vector $\hat{\vek x}_{0}\Sfs$, we now consider the
performance of GMRES applied to the shifted system.
\bthm\label{thm.GMRES-Ritz-performance}
	Let $\vek A\in\C^{n\times n}$ be diagonalizable with eigendecomposition \eqref{eqn.A-eig-decomp}.
	Let $\CK_{j}(\vek A,\vek r_{0})$ be the Krylov subspace generated by $j$ iterations of 
	unrestarted GMRES applied to \eqref{eqn.Axb}, and for an indexing set $\I\subsetneq\curl{1,\ldots, n}$
	let $\CY$ be an invariant subspace of $\vek A$
	spanned by $\curl{\vek f_{i}}_{i\in\I}$.  Let $\vek P_{\CY}$ be the orthogonal projection onto $\CY$, and let
	$\overline{\vek x}_{0}\Sfs$ be the result of the error projection 
	\be\nn
		\vek x\Sfs - \overline{\vek x}_{0}\Sfs = \prn{\vek I - \vek P_{\CY}}\prn{\vek x\Sfs - \hat{\vek x}_{0}\Sfs}.
	\ee
	
	If we apply Algorithm {\rm \ref{algorithm.sgmres}} to solve \eqref{eqn.model-problem-unprec} with
	no preconditioning and no restarting, then the residual $\hat{\vek r}_{\ell}\Sfs$ resulting from 
	$\ell$ iterations of GMRES applied to the shifted system with starting vector $\hat{\vek x}_{\ell}\Sfs$
	defined as in \eqref{eqn.error-projection} satisfies the bound
	\be\label{eqn.shift-GMRES-bound}
		\norm{\hat{\vek r}_{\ell}}\leq \norm{\overline{\vek r}_{\ell}\Sfs} + \delta	
	\ee
	where $\overline{\vek r}_{\ell}\Sfs$ is the residual resulting from applying $\ell$ iterations
	of GMRES to the shifted system with starting vector $\overline{\vek x}_{0}\Sfs$, and 
	$\delta = \sum_{i\in\I}(\lambda_{i}+\sigma)\overline{p}_{\ell}(\lambda_{i}+\sigma)\phi_{i}\vek f_{j}$,
	with $\overline{p}_{j}$ being the $j$th GMRES residual polynomial associated to the iteration 
	for starting vector $\overline{\vek x}_{0}\Sfs$ and 
	$\phi_{i}=\vek w_{i}^{\ast}\prn{\vek I - \vek P}\prn{\vek x\Sfs - \vek x_{0}\Sfs}$.
\ethm
\bproof
	The structure of this proof follows that in \cite{KMR.2001}, but it is also related to the results
	presented in \cite[Section 3.1, Case 1]{CN.GalProjAnalMRHS.1999} of CG with Lanczos-Galerkin projection
	applied to shifted systems.
	
	Because $\vek A$ is diagonalizable, we can decompose the errors with respect to
	$\hat{\vek x}\Sfs$ and $\overline{\vek x}\Sfs$ as 
	\be\nn
		\vek x\Sfs - \hat{\vek x}_{0}\Sfs = \sum_{i=1}^{n}\phi_{i}\vek f_{i}\mand\vek x\Sfs - \overline{\vek x}_{0}\Sfs = \sum_{\underset{i\not\in\I}{i=1}}^{n}\phi_{i}\vek f_{i},
	\ee
	which implies that 
	\be\nn
		\hat{\vek r}_{0}\Sfs = \sum_{i=1}^{n}\phi_{i}\prn{\lambda_{i}+\sigma}\vek f_{i}\mand\overline{\vek r}_{0}\Sfs = \sum_{\underset{i\not\in\I}{i=1}}^{n}\phi_{i}\prn{\lambda_{i}+\sigma}\vek f_{i}.
	\ee
	Because the GMRES residual polynomial $\hat{p}_{\ell}$ satisfies the minimization
	\be\nn
		\hat{p}_{\ell} = \argmin{\underset{p(0)=1}{p\in\Pi_{\ell}}}\norm{p(A+\sigma\vek I)\hat{r}_{0}}\mwith \Pi_{\ell} = \curl{p\,|\,\deg p \leq \ell},
	\ee
	we can write 
	\bea
		\norm{\hat{\vek r}_{\ell}} = \min_{\underset{p(0)=1}{p\in\Pi_{\ell}}}\norm{p(A+\sigma\vek I)\hat{r}_{0}} & \leq & 
		\norm{\overline{p}_{\ell}\prn{\vek A + \sigma}\hat{\vek r}_{0}} = \norm{\sum_{i=1}^{n}(\lambda_{i}+\sigma)\overline{p}_{\ell}(\lambda_{i} + \sigma)\phi_{i}\vek f_{i}}\nn\\
		&=&\norm{\sum_{\underset{i\not\in\I}{i=1}}^{n}(\lambda_{i}+\sigma)\overline{p}_{\ell}(\lambda_{i}+\sigma)\vek f_{i} + \sum_{i\in\I}(\lambda_{i}+\sigma)\overline{p}_{\ell}(\lambda_{i}+\sigma)\phi_{i}\vek f_{i}}\nn\\
		&\leq &\norm{\overline{\vek r}_{\ell}\Sfs}  + \underbrace{\norm{\sum_{i\in\I}(\lambda_{i}+\sigma)\overline{p}_{\ell}(\lambda_{i}+\sigma)\phi_{i}\vek f_{i}}}_{\delta}.\nn
	\eea
	From the definitions of $\vek f_{i_{1}}$ and $\vek w_{i_{2}}$, we know that 
	$\vek w_{i_{2}}^{\ast}\vek f_{i_{1}}=\delta_{i_{1},i_{2}}$.  Thus, from the definition of 
	$\vek x\Sfs - \hat{\vek x}_{0}\Sfs$ as well as its eigendecomposition, we have that \linebreak
	$\phi_{i}~=~\vek w_{i}^{\ast}\prn{\vek I-\vek P}\prn{\vek x\Sfs - {\vek x}_{0}\Sfs}$.
\eproof

Certainly, Theorem \ref{thm.GMRES-Ritz-performance} applies to any invariant subspace $\CY$.  However, the 
interesting case, which is considered in \cite{Chan1997,CN.GalProjAnalMRHS.1999,KMR.2001}, is when 
$\CY$ is such that the Krylov subspace $\CK_{j}(\vek A,\vek r_{0})$ contains a good approximation of it.
If $\CK_{j}(\vek A,\vek r_{0})$ actually contained $\CY$, then it is straightforward to show that 
$\phi_{i} = 0$ for all $i\in\I$, and thus $\delta = 0$.  We can then expect that if $\CY$ is well-approximated
in $\CK_{j}(\vek A,\vek r_{0})$, that $\delta$ would be non-zero but small.  In this case, the behavior
of GMRES applied to the shifted system with starting vector $\hat{\vek x}_{0}$ would mimic GMRES 
applied to that same system with with starting vector $\overline{\vek x}_{0}$, in which the iteration is 
orthogonal to $\CY$.  Unfortunately, this theory cannot be trivially extended to the case that the correction space
is not a Krylov subspace, as it relies on the polynomial approximation interpretation of GMRES.

Following from \cite{Saad1987}, we also can analyze the effectiveness of the direct projection by 
decomposing the residual. This analysis is developed in the general
framework setting presented in Section \ref{section.dir-proj-framework} 
and then interpreted for the individual methods.  
Here we use the notation that $\vek P(\cdot)$ denotes the orthogonal projector onto
the subspace specified in the argument.  
\bthm\label{thm.resid-decomp}
	Let the sequence
	of subspaces $\curl{\CS_{m}}$ be defined as in Section {\rm\ref{section.dir-proj-framework}},
	and additionally let
	\be\label{eqn.subset-containment-assumpt}
		\CT_{m}= \prn{\vek A + \sigma\vek I}\vek M^{-1}\CS_{m}.
	\ee
	If $\vek r_{0}\Sfs$ is the initial residual
	for the shifted system, and $\hat{\vek r}_{0}\Sfs$ is the residual produced by
	projecting $\vek r_{0}\Sfs$ according to \eqref{eqn.prec-Galerkin}, 
	then we have that 
	\be\label{eqn.shifted-resid-decomp}
		\hat{\vek r}_{0}\Sfs = \prn{\vek I - \vek P\prn{\CT_{m}}}\vek P\prn{\CT_{m+1}}\vek r_{0}\Sfs + \prn{\vek I - \vek P\prn{\CT_{m}}}\vek r_{0}\Sfs.
	\ee	 
\ethm
Note that in the unpreconditioned case, Theorem \ref{thm.resid-decomp} can be applied by taking $\vek M=\vek I$.
\bproof  
Using the property of projectors, we can decompose
\be\nn
	\vek r_{0}\Sfs = \vek P\prn{\CT_{m+1}}\vek r_{0}\Sfs + (\vek I - \vek P\prn{\CT_{m+1}})\vek r_{0}\Sfs.  
\ee
The minimum residual projection \eqref{eqn.prec-Galerkin} can be written,
\be
	\hat{\vek r}_{0}\Sfs = (\vek I - \vek P\prn{\CT_{m}})\vek P\prn{\CT_{m+1}}\vek r_{0}\Sfs + (\vek I - \vek P\prn{\CT_{m}})(\vek I - \vek P\prn{\CT_{m+1}})\vek r_{0}\Sfs.
\ee
From \eqref{eqn.subset-containment-assumpt} and the definition of 
$\CT_{m}$, we have that 
\be\nn
		\CT_{m}\subset \CT_{m+1}
\ee
which in turn yields the reverse containment of the orthogonal complements,
\be\nn
	{\CT_{m+1}}^{\perp}\subset\CT_{m}^{\perp}
\ee
and thus
\be\nn
	(\vek I - \vek P\prn{{\CT}_{m}})(\vek I - \vek P\prn{{\CT}_{m+1}})\vek r_{0}\Sfs = (\vek I - \vek P\prn{{\CT}_{m+1}})\vek r_{0}\Sfs 
\ee
This yields the result. \hfill 
\eproof

\bcor
	Let the same assumptions as in {\rm Lemma \ref{thm.resid-decomp}} hold.  Then we have the following bound on $\norm{\hat{\vek r}_{0}\Sfs}$,
	\be\label{eqn.resid-bound}
		\norm{\hat{\vek r}_{0}\Sfs} \leq \norm{\prn{\vek I - \vek P\prn{{\CT}_{m}}}\vek P\prn{\CT_{m+1}}\vek r_{0}\Sfs} + \norm{\prn{\vek I - \vek P\prn{\CT_{m+1}}}\vek r_{0}\Sfs}
	\ee
\ecor
\bproof 
We simply take the norm of both sides of \eqref{eqn.shifted-resid-decomp} and apply the triangle inequality.\hfill
\eproof

From \eqref{eqn.resid-bound}, we can see that the residual norm bound 
depends on both the effectiveness of the minimization
projection applied to the orthogonal projection of $\vek r_{0}\Sfs$ in 
$\CT_{m+1}$ and the size of the part of the residual which lies
in $\CT_{m+1}^{\perp}$.
As an aside, to connect this analysis back to the two proposed methods, 
we observe that in the case of the right-preconditioned
 shifted GMRES algorithm (Algorithm \ref{algorithm.sgmres}), we have 
 \begin{align}
 \CS_{m} &= \vek M^{-1}\CK_{m}(\vek A\vek M^{-1},\vek r_{0}),\mbox{\ \ }\CT_{m} = \vek A\vek M^{-1}\CK_{m}(\vek A\vek M^{-1},\vek r_{0}),\mbox{\ \ and,}\nn\\ \CT_{m} &= \vek A\vek M^{-1}\CK_{m}(\vek A\vek M^{-1},\vek r_{0}) + \sigma\vek M^{-1}\CK_{m}(\vek A\vek M^{-1},\vek r_{0}).
\end{align}  
In the case of the preconditioned rGMRES method for shifted systems 
(Algorithm \ref{algorithm.sgcrodr}), we have
{\footnotesize
\begin{align}
	\CS_{m} &=\vek M^{-1}\curl{\CU + \CK_{m}((\vek I - \vek P)\vek A\vek M^{-1}, \vek r_{0})},\mbox{\ \ }\CT_{m} = \CC + \vek A\vek M^{-1}\CK_{m}((\vek I - \vek P)\vek A\vek M^{-1}, \vek r_{0}), \mbox{\ \ and, }\nn\\ \CT_{m} &= \CC + \vek A\vek M^{-1}\CK_{m}((\vek I - \vek P)\vek A\vek M^{-1}, \vek r_{0}) + \sigma\vek M^{-1}\curl{\CU + \CK_{m}((\vek I - \vek P)\vek A\vek M^{-1},\vek r_{0})}
\end{align}}%

As a quick aside, we mention briefly that the matrix $\vek N_{m}\Sfs$ is connected 
to a generalized eigenvalue approximation problem associated to the 
computation of the harmonic Ritz values; see, e.g.,  
\cite{Morgan1998,Paige1995}.  This is elaborated upon in the tech report \cite{S.2014.2}.

\subsection{Cost of the algorithms}\label{section.cost}
For Algorithms \ref{algorithm.sgmres} and \ref{algorithm.sgcrodr}, we enumerate the additional 
per-cycle costs incurred by the proposed algorithms as they are built, respectively, on top of a cycle
of GMRES and a cycle of \RGMRES.  

Let $c_{old}$ denote the cost per iteration
of an existing method (here GMRES or \RGMRES) and $c_{new}$ the cost per iteration of the modified method
(here Algorithm \ref{algorithm.sgmres} or \ref{algorithm.sgcrodr}).  Here we don't specify how cost should be measured.
It could be by estimating, e.g., FLOPS, amount of data moved, actual timings of various operations, etc.  
In our subsequent calculations, though, we estimate costs in FLOPS.
In this setting, we have that the new methods
cost more per iteration, i.e., $c_{new} = c_{old} + d_{new}$.  In Tables \ref{table.sgmres-ops-dazu} and
\ref{table.sgcrodr-ops-dazu}, we list, respectively, the additional costs of each proposed algorithms, allowing us to estimate
$d_{new}$.  An important consideration which we don't treat here is the cost of applying the operator, which depends on
characteristics such as sparsity.  This  can dominate the cost per iteration.  In judging the effectiveness
of these methods, the benefit of iteration reduction is dictated by the matrix-vector product cost (which would also
include the cost of applying the preconditioner).

We can similarly define the number of iterations
required by both methods to solve all shifted systems, i.e., $j_{old}$ and $j_{new}$.  By assumption, 
the newer method should solve all shifted systems in fewer iterations, i.e., $j_{new} = j_{old} - a_{new}$.
Roughly speaking then, the total cost of each method can be estimated by $j_{old}\cdot c_{old}$ and $j_{new}\cdot c_{new}$.

\subsubsection{Comparison of Algorithm \ref{algorithm.sgmres} to GMRES}\label{section.sGMRES-cost}
Algorithm \ref{algorithm.sgmres} is built on top of GMRES.  In Table \ref{table.sgmres-ops-dazu},
we list all additional operations and information about their costs.
\input{sgmres-ops-dazu}
From this, we can estimate
the additional per cycle FLOP cost and then divide by $m$ to estimate $d_{new}^{\eqref{algorithm.sgmres}}$.  
If we simplify, we see that
\be\nn
	d_{new}^{\eqref{algorithm.sgmres}} = \frac{2}{3} (L+3) m^2+2 m (2 L+n+1)+6 L n+n.
\ee

\subsubsection{Comparison of Algorithm \ref{algorithm.sgcrodr} to \RGMRES}
Algorithm \ref{algorithm.sgcrodr} is built on top of
recycled GMRES.  We can compare costs of a cycle of each
algorithm by looking at the additional costs per cycle of Algorithm \ref{algorithm.sgcrodr}.  
There are also a few initial one-time overhead costs which must be taken
into account.  Thus in Table \ref{table.sgcrodr-ops-dazu} we show the additional per cycle
costs of Algorithm \ref{algorithm.sgcrodr}, and in Table \ref{table.sgcrodr-ops-dazu-onetime}
we show the additional one-time overhead costs.
\input{sgcrodr-ops-dazu}
We use Table \ref{table.sgcrodr-ops-dazu} to 
estimate $d_{new}^{\eqref{algorithm.sgcrodr}}$ but we must also take into account the onetime costs shown in 
Table \ref{table.sgcrodr-ops-dazu-onetime} by dividing those costs by 
the total number of iterations $j_{new}^{\eqref{algorithm.sgcrodr}}$.  
After simplifying we have
\bea
	d_{new}^{\eqref{algorithm.sgcrodr}} &=& (1 + \frac{5L}{3})m^{2} + (2 + 3 k + 3 L + 5 k L + 2 n)m \nn\\
					&&+ 1+4 k+3 k^2+4 L+6 k L+5 k^2 L+n+3 k n+6 L n\nn\\ 
					&&+ \frac{1}{m}(\frac{5 k^3 L}{3}+k^3+3 k^2 L+2 k^2 n+2 k^2+6 k L n+4 k L+k n+k) \nn\\
					&&+ \frac{k^3+2 k^2 n+\left(k^3+\frac{2 k^2}{3}\right) L+6 k L n+3 k L}{j_{new}}.\nn
\eea
\subsubsection{Estimating costs for specific examples}\label{eqn.cost-estimates-ex}
Now we can compare costs for a specific example.  For Algorithm \ref{algorithm.sgmres}, let $m=50$, 
$L=5$, and $n=10^{5}$.  Then we have 
$d_{new}^{\eqref{algorithm.sgmres}} \approx 1.3\times 10^{6}$.  For Algorithm
\ref{algorithm.sgcrodr}, let us store a small recycled subspace but use the same amount of storage, i.e., 
$m=40$ and $k=5$.  This yields \linebreak $d_{new}^{\eqref{algorithm.sgcrodr}}\approx 1.3\times 10^{6}  + \frac{2.0\times 10^{6}}{j_{new}^{\eqref{algorithm.sgcrodr}}}$.

Admittedly, this is a bit unwieldy and has many parameters.  However, if we make an additional assumption on
how Algorithm \ref{algorithm.sgcrodr} is called, we can simplify the associated cost calculation.  Let us assume
that $k=\frac{1}{2}m$, i.e., that we maintain a recycled subspace half the size of the associated Krylov subspace
dimension. 
Then we see that we can simplify
\bea
	d_{new}^{\eqref{algorithm.sgcrodr}} &=& m^3 \left(\frac{L}{8 j_{new}^{\eqref{algorithm.sgcrodr}}}+\frac{1}{8 j_{new}^{\eqref{algorithm.sgcrodr}}}\right)+ m^2 \left(\frac{L}{6 j_{new}^{\eqref{algorithm.sgcrodr}}}+\frac{45 L}{8}+\frac{n}{2 j_{new}^{\eqref{algorithm.sgcrodr}}}+\frac{27}{8}\right)\nn\\ 
					&& + m \left(\frac{3 L n}{j_{new}}+\frac{3 L}{2 j_{new}^{\eqref{algorithm.sgcrodr}}}+\frac{27 L}{4}+4 n+\frac{9}{2}\right)+ \frac{3}{2}+6 L+\frac{3 n}{2}+9 L n.\nn
\eea
Let us assume for Algorithm \ref{algorithm.sgmres} that we have the same values as before.  To have 
approximately equivalent storage for Algorithm \ref{algorithm.sgcrodr}, we set $m=25$,  and we have  
$d_{new}^{\eqref{algorithm.sgcrodr}} \approx 1.5\times 10^{6}  + \frac{6.9\times 10^{6}}{j_{new}^{\eqref{algorithm.sgcrodr}}}$ FLOPS.

\input{cost_calc.tex}

In Figure \ref{figure.cost_calc}, we study the growth in estimated FLOP costs when all but one parameter are held fixed.
For srGMRES, we again assume that $k=\frac{1}{2}m$ and that 
the total number of iterations needed for Algorithm \ref{algorithm.sgcrodr} to converge is
\be\nn
	j_{new}^{\eqref{algorithm.sgcrodr}} = \frac{n}{{10^{\frac{4}{9}+\frac{2m}{9m+9}}\log_{10}n}}.
\ee
This is somewhat arbitrary, but it qualitatively matches experimental
observations.  This formula is derived so that for the case that $n=10^{6}$ and for $m\rightarrow\infty$ we have
that $j_{new}^{\eqref{algorithm.sgcrodr}}\rightarrow 100$ and the convergence is monotonically decreasing and relatively
fast.
It is necessary to have some assumption on the value of $j_{new}^{\eqref{algorithm.sgcrodr}}$
since Algorithm \ref{algorithm.sgcrodr} has some overhead costs which need to be amortized over the total number 
of iterations.   In the three graphs shown; we vary, respectively, number of shifts ($L$), problem dimension ($n$), and 
cycle length ($m$) with everything else being held constant.  For the experiments in which $L$ is held constant, we chose
$L=5$. Similarly, we chose $n=10^{7}$ and $m=100$ in the cases that these parameters were held constant.

We conclude by noting that we consider only one type of costs in this section.  In reality, these methods
also incur storage costs and data movement costs which are nontrivial for large-scale problems
which must be considered.  Furthermore, absent preconditioning, it is clear from the cost calculations
that in the case of non-Hermitian shifted systems of the form treated in \cite{Kilmer.deSturler.tomography.2006}
that the method considered in that paper would be much cheaper than Algorithm \ref{algorithm.sgcrodr}, 
and absent preconditioning,
for general non-Hermitian
shifted linear systems satisfying the conditions in \cite{Frommer1998} (e.g., collinear residuals), 
that method would outperform Algorithm \ref{algorithm.sgmres}.
Lastly, both Algorithms \ref{algorithm.sgmres} and \ref{algorithm.sgcrodr} can be used with flexible
preconditioners with \emph{no additional computational or storage costs}.

\section{Numerical Results}\label{section.num-results}
We performed a series of numerical experiments to demonstrate the effectiveness of our algorithms as
well as to compare performance (as measured in both matrix-vector product counts and CPU timings) with
other algorithms.  
All tests were performed in Matlab R2014b (8.4.0.150421) 64-bit running on a Mac Pro workstation with
two $2.26$ GHz Quad-Core Intel Xeon processors and 12 GB 1066 MHz DDR3 main memory.
For these tests, we use \emph{two sets of QCD matrices} downloaded from the University of 
Florida Sparse Matrix Library \cite{DH.2011}.  One set of matrices is a collection of 
seven $3072\times 3072$ complex matrices and the other is a collection of seven 
$49152\times 49152$ complex matrices.  
For each matrix $\vek D$ from the collection, there exists some critical value $\kappa_{c}$ 
such that for  $\frac{1}{\kappa_{c}}< \frac{1}{\kappa}<\infty$, the matrix 
$\vek A =  \frac{1}{\kappa}\vek I - \vek D$ is real-positive. 
For each $\vek D$, we took $\vek A = \prn{\frac{1}{\kappa_{c}} + 10^{-3}}\vek I - \vek D$
as our base matrix.
In our experiments then, each set is taken as the sequence $\curl{\vek A_{i}}$ and we solve a family
of the form \eqref{eqn.one-shifted-family-seq}.  As described in \cite{DH.2011}, the matrices
$\vek D$ are discretizations of the Dirac operator used in numerical
			simulation of quark behavior at different physical temperatures.
			We note that larger real shifts of $\vek A_{i}$ yields better conditioned matrices for all $i$.
For all experiments, we chose the right-hand side
$\vek b_{1} = \vek 1$, the vector of ones and set $\vek b_{i} = \vek b_{i-1} + \vek d_{i}$ where $\vek d_{i}$
is chosen randomly such that $\norm{\vek d_{i}}=1e-1$.
The requested relative residual tolerance was
$\varepsilon =10^{-8}$.  All augmentation was with harmonic Ritz vectors.
For all experiments, we preconditioned with an incomplete
LU-factorization (ILU) for the system with the smallest shift,
constructed using the Matlab function \verb|ilu()| called with the 
default Matlab settings.  We comment that the usage of ILU was a matter of 
convenience and effectiveness for these sample problems.  Its usage is meant to
demonstrate proof-of-concept rather than as advocating the usage of ILU 
for large-scale QCD problems.

We also comment about methods which we have omitted from testing;
the shifted restarted GMRES method \cite{Frommer1998}, the shifted GMRES-DR method \cite{Darnell2008}, and the 
recursive \RGMRES\ method for shifted systems proposed in 
\cite{SSX.2014}.  We have omitted these methods from the tests as they
do not admit general preconditioning. As such, they require substantially
more iterations in many experiments.  However, with the methods in \cite{Darnell2008,Frommer1998}, 
there would be some number of shifts for which this method
would be superior to those presented in this paper, as
cost of recursion in our methods, even with preconditioning,
would be greater than simply solving the unpreconditioned problems simultaneously with their shifted GMRES method \cite{Frommer1998}.

Since these experiments involve solving shifted systems
with shifts of varying magnitudes, it is useful to know information about the norms of our test matrices.
Therefore, we provide both the one- and two-norms for these matrices 
(computed respectively with the Matlab functions 
{\tt norm(}$\cdot${\tt , 1)} and {\tt svds(}$\cdot${\tt , 1)}).
The $1$-norms of these matrices all lie in the interval $\prn{28,31}$, and their $2$-norms lie in the interval
$\prn{11,14}$.

In our first experiment, we tested Algorithm \ref{algorithm.sgcrodr} with \emph{the smaller set of matrices}
 for various recycle space dimension sizes and restart cycle
lengths.  We solve for shifts $\sigma\in\curl{.01, .02, .03, 1, 2, 3}$.
We calculated total
required matrix-vector products.  
We see in Table~\ref{table.qcdSmallmkTable} that for these particular QCD matrices, good results can be achieved
for a small recycled subspace dimension as long as the cycle length is sufficiently
long.
\input{qcdSmallmkTable.tex}

For the remaining tests, 
we use the \emph{larger set of QCD matrices}.
In Table \ref{table.qcdLargeTimingComparison}
we compare time and matrix-vector product counts.  
For Algorithm \ref{algorithm.sgcrodr}, we chose cycle-length/recycle subspace dimension pair $(m,k)=(80,10)$ 
and use this pair for all experiments with Algorithm \ref{algorithm.sgcrodr} except for the one shown in 
Figure \ref{figure.qcdBigRecycleDimsPlot}.  
Parameters for Algorithm \ref{algorithm.sgmres} and other tested methods were chosen in order to have 
the same per-cycle storage cost of $3k+2m=190$ vectors \footnote[4]{for storing $\vek V_{m}$, $\vek Z_{m}$,$\vek U$,
$\vek C$, and $\vek Z_{\CU}$}.
For each family of linear systems, the experiment was performed ten times and the 
average time over these ten runs was taken as the run time.
We solved for a larger number of shifts of varying magnitudes, 
 $$\sigma\in\curl{.001, .002, .003, .04,.05,.06,.07, .8, .9, 1, 1.1, 10, 11, 12}.$$
We compared four methods (Algorithm \ref{algorithm.sgmres}, 
Algorithm \ref{algorithm.sgcrodr}, sequentially
applied GMRES and sequentially applied \RGMRES).
We see that for this problem with these shifts, both proposed algorithms outperform
the sequential applications of GMRES and \RGMRES\ both in terms of matrix-vector product counts
and run times.  In this case, the sGMRES algorithm is superior in time to 
srGMRES but not in terms of matrix-vector products, which demonstrates the difference in overhead costs.  
\input{qcdLargeTimingComparison.tex}

In Figure \ref{figure.qcdBigRecycleDimsPlot}, for a total fixed augmented subspace
dimension of $100$, we investigate how many matrix vector products are required
to solve the same sequence of problems with the same shifts as in the 
previous experiment for different values of $(m,k)$ 
such that $m+k=100$ where $m$ is the dimension
of the projected Krylov subspace and $k$ is the dimension of the recycled 
subspace.  With this we demonstrate a reduction in iterations as we allow 
more information to be retained in the subspace.
\input{qcdBigRecycleDimsPlot}

In Table \ref{table.qcdLargeShiftSizeMatvecCounts}, we study
matrix-vector product counts for different methods for shifts of 
varying magnitudes.  For each shift, we solve just two systems, the base system
and one shifted system.  Thus we can see how many additional matrix-vector products
are required for shifts of different magnitudes.
What we see is that for this set of matrices, overall performance
does not depend on shift magnitude.  For larger shifts, we see that 
Algorithm~\ref{algorithm.sgcrodr} and sequentially applied rGMRES are comparable
\textbf{when there is only one shift}.  For the QCD matrices, larger real shifts produce better 
conditioned problems and Table \ref{table.qcdLargeShiftSizeMatvecCounts} illustrates the trade-off
between better conditioning and reduced effectiveness of the proposed algorithm for larger shifts.
We hypothesize that the smallest values that are attained in the middle of the table are the result of 
Algorithms \ref{algorithm.sgmres} and \ref{algorithm.sgcrodr} still being effective for $\CO(1)$ shifts 
where we also see improved conditioning of the shifted systems.
\input{qcdLargeShiftSizeMatvecCounts.tex}

However, we have seen for larger numbers of shifts that 
Algorithm~\ref{algorithm.sgcrodr} can exhibit superior performance.  This raises the 
question, what are the marginal costs of solving each additional linear system
for \RGMRES\ and shifted \RGMRES, i.e., how many more matrix-vector products 
does each new shifted system require?  
This is investigated in Figure \ref{figure.qcdLargeMarginalCost}.
For two sets of twenty shifts, we calculated the marginal 
cost of solving each additional shifted system using Algorithm \ref{algorithm.sgcrodr}
as compared to \RGMRES.  
In Figure \ref{figure.qcdLargeMarginalCost} the first set of shifts (left-hand figure) were evenly space points
from the interval $\brac{0,1}$, and the second set of shifts (right-hand figure)
were evenly spaced points
from the larger interval $\brac{1,10}$.
In Figure \ref{figure.qcdLargeMarginalCost},
\input{qcdLargeMarginalCost}
we see that for the smaller interval, the cost of each new shifted system drops
for both algorithms but that Algorithm \ref{algorithm.sgcrodr} has the lower marginal
cost per shift.  For the larger set of shifts, we see that the marginal costs for both
algorithms actually increases for each new shift.  However, 
the marginal cost of each 
new shifted system for Algorithm \ref{algorithm.sgcrodr} 
becomes more stable (it levels off).
For sequentially applied \RGMRES, the marginal costs increases steadily for
all twenty shifts.  

In Figure \ref{figure.shiftMagnitudeConvergenceComparison}, 
we show the residual histories for systems solved using 
Algorithm~\ref{algorithm.sgcrodr} for shifts of various 
magnitudes, 
$$\sigma\in\curl{10^{-3},10^{-2},10^{-1},1,10^{1},10^{2},10^{3}}.$$
When viewing Figure 
\ref{figure.shiftMagnitudeConvergenceComparison}, we see (in this example) 
that the amount of improvement
\input{shiftMagnitudeConvergenceComparison}
for the shifted residuals is somewhat predicted by the shift magnitude,
though we again observe that the better conditioning of the systems with larger shifts seems to lead to 
more rapid convergence but at the expense of reduced effectiveness of the Lanczos-Galerkin projection.

Omitted here is a study of the
eigendecomposition of the residuals, which yielded
no discernible damping of certain eigenmodes
or other interesting observable phenomena after the projection
of the shifted residuals in our experiments.  Such experiments were to investigate questions of 
the convergence rates observed in Figure \ref{figure.shiftMagnitudeConvergenceComparison}.

\section{Conclusions}\label{section.conclusions}
We have presented two new methods for solving a family or a sequence of 
families of shifted linear systems with general preconditioning.  These methods
are derived from a general framework, which we also developed in this paper.
  These methods use subspaces generated during the minimum residual iteration 
  of the
base system to perform the projections for the shifted systems.  This 
technique
is fully compatible with right preconditioning, requiring only some additional 
storage.
The strength of methods derived from this framework is that 
preconditioned methods for shifted systems
easily can be built on top of existing minimum residual projection algorithms (and existing
codes) with only minor modifications.  We developed 
two algorithms: shifted GMRES
and shifted \RGMRES.  We demonstrated with numerical experiments that both 
methods can perform competitively.  

Finally, we note that our framework
is fully compatible with flexible and inexact Krylov subspace methods.
As this work all follows from \cite{CN.GalProjAnalMRHS.1999,Chan1997,KMR.2001,Saad1987}, it is also
clear that the method is also applicable to the case that we are solving \eqref{eqn.one-shifted-family-seq}
but with right-hand sides $\vek b_{i,\ell}$ which vary both with respect to coefficient 
matrix $\vek A_{i}$ and shift $\sigma_{\ell}$.

\section*{Acknowledgments}
 The author would like to thank Michael Parks who, while reviewing the author's dissertation, made a 
 comment which inspired this work.  The author would also like to thank Valeria Simoncini for insightful questions and comments during the author's visit to Bologna and Daniel Szyld for constructive comments. The author further thanks both 
 reviewers for offering comments and criticisms which led to great improvement in the presentation and
 completeness of this work.

\bibliographystyle{siam}
\bibliography{master}
\end{document}

%% file: macros.tex
\def\cal{\mathcal}
\def\norm#1{\left\|#1\right\|}

\def\C{\mathbb{C}}
\def\N{\mathbb{N}}

\def\I{\mathbb{I}}
\def\Cn{\C^n}

\def\Cnn{\C^{n\times n}}

\def\BA{{\bf A}}  
  
  \def\CC{{\cal C}}

  \def\CK{{\cal K}}

  \def\CO{{\cal O}}
  \def\CP{{\cal P}}

  \def\CS{{\cal S}}
  \def\CT{{\cal T}}
  \def\CU{{\cal U}}

  \def\CY{{\cal Y}}

\def\BAs{\BA{\kern-1.5pt}}

\def\CPs{\CP{\kern-0.8pt}}

\mathcode`\@="8000
{\catcode`\@=\active \gdef@{\mkern1mu}}

\def\mydate{\number\day\ {\ifcase\month \or January\or February\or
              March\or April\or May\or June\or July\or August\or
              September\or October\or November\or December\fi}
\number\year}
\def\vek#1{\mathbf{#1}}

\newcommand{\argmin}[1]{\underset{#1}{\text{{\rm argmin}}}}

\newcommand{\prn}[1]{\left(#1\right)}

\newcommand{\brac}[1]{\left[#1\right]}
\newcommand{\curl}[1]{\left\{#1\right\}}
\newcommand{\ab}[1]{\left|#1\right|}

\newcommand\restr[2]{{
  \left.\kern-\nulldelimiterspace 
  #1 
  \vphantom{\big|} 
  \right|_{#2} 
  }}

\def\be{\begin{equation}}
\def\ee{\end{equation}}
\def\bea{\begin{eqnarray}}
\def\eea{\end{eqnarray}}
\def\nn{\nonumber}
\def\mand{\mbox{\ \ \ and\ \ \ }}

\def\mwith{\mbox{\ \ \ with\ \ \ }}

\def\bbmat{\begin{bmatrix}}
\def\ebmat{\end{bmatrix}}

\def\balg{\begin{algorithm}}
\def\ealg{\end{algorithm}}
\def\bthm{\begin{theorem}}
\def\ethm{\end{theorem}}
\def\blem{\begin{lemma}}
\def\elem{\end{lemma}}
\def\bprop{\begin{proposition}}
\def\eprop{\end{proposition}}
\def\bcor{\begin{corollary}}
\def\ecor{\end{corollary}}
\def\bdefin{\begin{definition}}
\def\edefin{\end{definition}}
\def\bc{\begin{cases}}
\def\ec{\end{cases}}
\def\bproof{\par\addvspace{1ex} \indent\textit{Proof}.\ \ }
\def\eproof{\hfill\cvd\linebreak\indent}
\newtheorem{exple}{Example}
\def\bex{\begin{exple}}
\def\eex{\end{exple}}
\newtheorem{assumption}{Assumption}
\def\bass{\begin{assumption}}
\def\eass{\end{assumption}}
\def\diag{{\rm diag}}

\def\cvd{~\vbox{\hrule\hbox{%
  \vrule height1.3ex\hskip0.8ex\vrule}\hrule } }

%% file: sgmres-algorithm.tex
\begin{algorithm}[hb!]
\caption{Right preconditioned shifted GMRES ({\tt sGMRES()})}
\label{algorithm.sgmres}
\SetKwInOut{Input}{Input}\SetKwInOut{Output}{Output}\SetKwComment{Comment}{}{}
\Input{$\vek A\in\Cnn$; $\vek b\in\Cn$; $\curl{\sigma_{\ell}}_{\ell=1}^{L}\subset\C$; Initial Approximations $\curl{\vek x\SfAlg{\sigma_{\ell}}}_{\ell=1}^{L}$; $\varepsilon > 0$; Cycle length $m\in\N$}
\Output{$\curl{\vek x\SfAlg{\sigma_{\ell}}}_{\ell=1}^{L}$ such that $\norm{\vek r\SfAlg{\sigma_{\ell}}}/\norm{\vek r_{0}\SfAlg{\sigma_{\ell}}}\leq\varepsilon$ for all $\ell$}
\For{$\ell = 1\ldots L$}{
	$\vek r\SfAlg{\sigma_{\ell}} = \vek b - (\vek A + \sigma_{\ell}\vek I)\vek x\SfAlg{\sigma_{\ell}}$
}
$\gamma_{1} = \norm{\vek r\SfAlg{\sigma_{1}} }$\\
\If{$L > 1$}{
	\While{$\norm{\vek r\SfAlg{\sigma_{1}}}/\gamma_{1}>\varepsilon$}{
		Compute and overwrite $\vek x\SfAlg{\sigma_{1}}$, $\vek r\SfAlg{\sigma_{1}}$, $\vek V_{m+1}$, $\vek Z_{m}$, $\overline{\vek H}_{m}$ by calling {\tt GMRES()} for $\vek A + \sigma_{1}\vek I$, $\vek M$, $\vek b$, $\vek x\SfAlg{\sigma_{1}}$, and $m$\\
		Compute and overwrite $\overline{\vek H}_{m}^{\ast}\overline{\vek H}_{m} $, 
 $\overline{\vek H}_{m}^{\ast}\vek V_{m+1}^{\ast}\vek Z_{m}$, and $\vek Z_{m}^{\ast}\vek Z_{m}$\label{line.3onetimeOps}\\
 		\For{$\ell=2\ldots L$}{
 			$\vek N \leftarrow \overline{\vek H}_{m}^{\ast}\overline{\vek H}_{m}+ \sigma\overline{\vek H}_{m}^{\ast}\vek V_{m+1}^{\ast}\vek Z_{m} + \overline{\sigma}\vek Z^{\ast}\vek V_{m+1}\overline{\vek H}_{m} + \ab{\sigma}^{2}\vek Z_{m}\vek Z_{m}$\label{line.Nsum}\\
 			$\vek y \leftarrow \vek N^{-1}\brac{\prn{\vek V_{m+1}\overline{\vek H}_{m} + \sigma \vek Z_{m}}}^{\ast}\vek r_{0}\Sfs$\label{line.Ninv-apply}\\
 			$\vek x\SfAlg{\sigma_{\ell}} \leftarrow \vek x_{0}\SfAlg{\sigma_{\ell}} + \vek Z_{m}\vek y$\label{line.x-update}\\
 			$\vek r\SfAlg{\sigma_{\ell}} \leftarrow \vek r_{0}\SfAlg{\sigma_{\ell}} - \prn{\vek V_{m+1}\overline{\vek H}_{m} + \sigma \vek Z_{m}}\vek y$\label{line.r-update}
 		}
 		For all $\ell = 2,\ldots L$ compute and overwrite $\vek x\SfAlg{\sigma_{\ell}}$ by recursively calling {\tt sGMRES()} for $\vek A$, $\vek b$, $\vek M$, $\curl{\sigma_{\ell}}_{\ell=2}^{L}$, $\curl{\vek x\SfAlg{\sigma_{\ell}}}_{\ell=2}^{L}$, $\varepsilon$, and $m$
	}
}\Else{
	\While{$\norm{\vek r\SfAlg{\sigma_{1}}}/\gamma_{1}>\varepsilon$}{
		Compute and overwrite $\vek x\SfAlg{\sigma_{1}}$, $\vek r\SfAlg{\sigma_{1}}$ by calling {\tt GMRES()} for $\vek A + \sigma_{1}\vek I$, $\vek M$, $\vek b$, $\vek x\SfAlg{\sigma_{1}}$, and $m$
	}
}
\end{algorithm}

%% file: sgcrodr-algorithm.tex
\begin{algorithm}[hb!]
\caption{Right preconditioned shifted \RGMRES\ ({\tt srGMRES()})}
\label{algorithm.sgcrodr}
\SetKwInOut{Input}{Input}\SetKwInOut{Output}{Output}\SetKwComment{Comment}{}{}
\Input{$\vek A\in\Cnn$; $\vek b\in\Cn$; $\curl{\sigma_{\ell}}_{\ell=1}^{L}\subset\C$; Initial Approximations $\curl{\vek x\SfAlg{\sigma_{\ell}}}_{\ell=1}^{L}$; $\vek U\in\C^{n\times k}$; $\varepsilon > 0$; Cycle length $m\in\N$}
\Output{$\curl{\vek x\SfAlg{\sigma_{\ell}}}_{\ell=1}^{L}$ such that $\norm{\vek r\SfAlg{\sigma_{\ell}}}/\norm{\vek r_{0}\SfAlg{\sigma_{\ell}}}\leq\varepsilon$ for all $\ell$}
\For{$\ell = 1\ldots L$}{
	$\vek r\SfAlg{\sigma_{\ell}} = \vek b - (\vek A + \sigma_{\ell}\vek I)\vek x\SfAlg{\sigma_{\ell}}$
}
$\gamma_{1} = \norm{\vek r\SfAlg{\sigma_{1}} }$\\
$\vek Z_{\CU} = \vek M^{-1}\vek U$\\
$\vek C = (\vek A + \sigma_{1}\vek I)\vek Z_{\CU}$\\
Compute QR-factorization $\vek Q\vek R = \vek C$\\
$\vek C \leftarrow \vek Q$, $\vek U \leftarrow \vek U\vek R^{-1}$, $\vek Z_{\CU} \leftarrow \vek Z_{\CU}\vek R^{-1}$\label{line.Zu-update}\\
$\vek x\SfAlg{\sigma_{1}} \leftarrow \vek x\SfAlg{\sigma_{1}}  + \vek U\vek C^{\ast}\vek r\SfAlg{\sigma_{1}} $ and $\vek r\SfAlg{\sigma_{1}} \leftarrow \vek r\SfAlg{\sigma_{1}}  - \vek C\vek C^{\ast}\vek r\SfAlg{\sigma_{1}} $ \\
Compute $\vek C^{\ast}\vek Z_{\CU} $ and $\vek Z_{\CU}^{\ast}\vek Z_{\CU}$\label{line.CTZu-ZuTZu}\\
\For{$\ell = 2\ldots L$}{
	\Comment{\%\%\%\%\% Shifted System Initial Projections \%\%\%\%\%}
	$\vek N \leftarrow \vek I + \sigma\vek C^{\ast}\vek Z_{\CU} + \overline{\sigma}\vek Z_{\CU}^{\ast}\vek C + \ab{\sigma}^{2}\vek Z_{\CU}^{\ast}\vek Z_{\CU}$\label{line.initial-N-sum}\\
	$\vek y \leftarrow \vek N^{-1}\prn{\vek C + \sigma\vek Z_{\CU}}^{\ast}\vek r\SfAlg{\sigma_{\ell}}$\label{line.initN-invert-update}\\
	$\vek x\SfAlg{\sigma_{\ell}} \leftarrow \vek x\SfAlg{\sigma_{\ell}} + \vek U\vek y$\label{line.initx-update}\\ 
	$\vek r\SfAlg{\sigma_{\ell}} \leftarrow \vek r\SfAlg{\sigma_{\ell}} - \prn{\vek C + \sigma\vek Z_{\CU}}\vek y$\label{line.initr-update}
}
\If{$L > 1$}{
	\While{$\norm{\vek r\SfAlg{\sigma_{1}}}/\gamma_{1}>\varepsilon$}{
		Compute and overwrite $\vek x\SfAlg{\sigma_{1}}$, $\vek r\SfAlg{\sigma_{1}}$, $\vek V_{m+1}$, $\vek Z_{m}$, $\overline{\vek H}_{m}$, $\vek B_{m}$ by calling {\tt rGMRES()} for $\vek A + \sigma_{1}\vek I$, $\vek M$, $\vek b$, $\vek x\SfAlg{\sigma_{1}}$, $\vek U$, $\vek C$, and $m$\\
		Compute and overwrite $\overline{\vek G}_{m}$, $\overline{\vek G}_{m}^{\ast}\overline{\vek G}_{m}$, $\vek Z_{\CU}^{\ast}\vek Z_{\CU}$, $\vek Z_{\CU}^{\ast}\vek Z_{m}$, $\vek Z_{m}^{\ast}\vek Z_{m}$, $\vek C^{\ast}\vek Z_{\CU}$, $\vek C^{\ast}\vek Z_{m}$, $\vek V_{m+1}^{\ast}\vek Z_{\CU}$, $\vek V_{m+1}^{\ast}\vek Z_{m}$\label{line.build-proj-mats}\\
 		\For{$\ell=2\ldots L$}{
 			\Comment{\%\%\%\%\% Shifted System Projections \%\%\%\%\%}
 			$\vek N \leftarrow \overline{\vek G}_{m}^{\ast}\overline{\vek G}_{m} + \ab{\sigma}^{2}\begin{bmatrix}\vek Z_{\CU}^{\ast}\vek Z_{\CU} & \vek Z_{\CU}^{\ast}\vek Z_{m}\\ \vek Z_{m}^{\ast}\vek Z_{\CU} & \vek Z_{m}^{\ast}\vek Z_{m}\end{bmatrix}	 
	+ \sigma\overline{\vek G}_{m}^{\ast}\begin{bmatrix}\vek C^{\ast}\vek Z_{\CU} & \vek C^{\ast}\vek Z_{m} \\ \vek V_{m+1}^{\ast}\vek Z_{\CU} & \vek V_{m+1}^{\ast}\vek Z_{m}\end{bmatrix} 
		+ \overline{\sigma}\begin{bmatrix}\vek Z_{\CU}^{\ast}\vek C & \vek Z_{\CU}^{\ast}\vek V_{m+1} \\ \vek Z_{m}^{\ast}\vek C & \vek Z_{m}^{\ast}\vek V_{m+1}\end{bmatrix}\overline{\vek G}_{m}$\label{line.build-N}\\
 			$\vek y \leftarrow \vek N^{-1}\curl{\begin{bmatrix}\vek C & \vek V_{m+1}\end{bmatrix}\overline{\vek G}_{m} + \sigma \begin{bmatrix}\vek Z_{\CU} & \vek Z_{m}\end{bmatrix}}^{\ast}\vek r\SfAlg{\sigma_{\ell}}$\label{line.N-inv}\\
 			$\vek x\SfAlg{\sigma_{\ell}} \leftarrow \vek x_{0}\SfAlg{\sigma_{\ell}} + \begin{bmatrix}\vek Z_{\CU} & \vek Z_{m}\end{bmatrix}\vek y$\label{line.update-x}\\
 			$\vek r\SfAlg{\sigma_{\ell}} \leftarrow \vek r_{0}\SfAlg{\sigma_{\ell}} - \curl{\begin{bmatrix}\vek C & \vek V_{m+1}\end{bmatrix}\overline{\vek G}_{m} + \sigma \begin{bmatrix}\vek Z_{\CU} & \vek Z_{m}\end{bmatrix}}\vek y$\label{line.update-r}
 		}
 		Compute updated $\vek U$, $\vek Z_{\CU}$, and $\vek C$
	}
	For all $\ell = 2,\ldots L$ compute and overwrite $\vek x\SfAlg{\sigma_{\ell}}$ by recursively calling {\tt srGMRES()} for $\vek A$, $\vek b$, $\vek M$, $\curl{\sigma_{\ell}}_{\ell=2}^{L}$, $\curl{\vek x\SfAlg{\sigma_{\ell}}}_{\ell=2}^{L}$, $\vek U$, $\varepsilon$, and $m$
}\Else{
	\While{$\norm{\vek r\SfAlg{\sigma_{1}}}/\gamma_{1}>\varepsilon$}{
		Compute and overwrite $\vek x\SfAlg{\sigma_{1}}$, $\vek r\SfAlg{\sigma_{1}}$ by calling {\tt rGMRES()} for $\vek A + \sigma_{1}\vek I$, $\vek M$, $\vek b$, $\vek x\SfAlg{\sigma_{1}}$, and $m$
	}
}
\end{algorithm}

%% file: sgmres-ops-dazu.tex
\begin{table}[htb]
\caption{Cost per cycle of extra calculations performed in Algorithm \ref{algorithm.sgmres} when compared to GMRES.}
\label{table.sgmres-ops-dazu}
\begin{center}
\begin{tabular}{r|c|c|c}
Operations & Alg. Line & FLOPS in $\CO(\cdot)$ & $\times$ per cycle \\
\hline
$\overline{\vek H}_{m}^{\ast}\overline{\vek H}_{m}$ & \ref{line.3onetimeOps} & $m^{3}+m^{2}$ & $1$\\
$\vek V_{m+1}^{\ast}\vek Z_{m}$ & \ref{line.3onetimeOps} & $n\prn{m^{2}+m}$ & $1$\\
$\overline{\vek H}_{m}^{\ast}\prn{\vek V_{m+1}^{\ast}\vek Z_{m}}$ & \ref{line.3onetimeOps} & $m^{3}+m^{2}$ & $1$\\
$\vek Z_{m}^{\ast}\vek Z_{m}$ & \ref{line.3onetimeOps} & $nm^{2}$ & $1$\\
Sum of 4 $m\times m$ matrices & \ref{line.Nsum} & $3m^{2}$ & $L$\\
$\brac{\prn{\vek V_{m+1}\overline{\vek H}_{m} + \sigma \vek Z_{m}}}^{\ast}\vek r_{0}\Sfs$ & \ref{line.Ninv-apply} & 2nm & $L$\\
Apply $\vek N^{-1}$ & \ref{line.Ninv-apply} & $\frac{2}{3}m^{3} + m^{2}$ & $L$\\
$\vek x\SfAlg{\sigma_{\ell}} \leftarrow \vek x_{0}\SfAlg{\sigma_{\ell}} + \vek Z_{m}\vek y$ & \ref{line.x-update} & 2mn & $L$\\
$\vek r\SfAlg{\sigma_{\ell}} \leftarrow \vek r_{0}\SfAlg{\sigma_{\ell}} - \prn{\vek V_{m+1}\overline{\vek H}_{m} + \sigma \vek Z_{m}}\vek y$ & \ref{line.r-update} & 2mn & $L$\\
\end{tabular}
\end{center}
\end{table}

%% file: sgcrodr-ops-dazu.tex
\begin{table}[htb]
\caption{Cost per cycle of extra calculations performed in Algorithm \ref{algorithm.sgcrodr} when compared to \RGMRES.}
\label{table.sgcrodr-ops-dazu}
\begin{center}
\begin{tabular}{r|c|c|c}
Operations & Alg. Line & FLOPS in $\CO(\cdot)$ & $\times$ per cycle \\
\hline
$\overline{\vek G}_{m}^{\ast}\overline{\vek G}_{m}$&\ref{line.build-proj-mats}& $(m+k+1)^{2}(m+k)$&1\\
$\vek Z_{\CU}^{\ast}\vek Z_{\CU}$&\ref{line.build-proj-mats}&$k^2$n&1\\
$\vek Z_{\CU}^{\ast}\vek Z_{m}$&\ref{line.build-proj-mats}&$knm$&1\\
$\vek Z_{m}^{\ast}\vek Z_{m}$&\ref{line.build-proj-mats}&$m^2n$&1\\
$\vek C^{\ast}\vek Z_{\CU}$&\ref{line.build-proj-mats}&$k^2n$&1\\
$\vek C^{\ast}\vek Z_{m}$&\ref{line.build-proj-mats}&$knm$&1\\
$\vek V_{m+1}^{\ast}\vek Z_{\CU}$&\ref{line.build-proj-mats}&$kn(m+1)$&1\\
$\vek V_{m+1}^{\ast}\vek Z_{m}$&\ref{line.build-proj-mats}&$nm(m+1)$&1\\
$\begin{bmatrix}\vek C & \vek V_{m+1}\end{bmatrix}\overline{\vek G}_{m} + \sigma \begin{bmatrix}\vek Z_{\CU} & \vek Z_{m}\end{bmatrix}$&\ref{line.N-inv}&$(m+k+1)^{2}(m+k)$&L\\
Sum of $4$ matrices &\ref{line.N-inv}&3(m+k)&L\\
\footnotesize$\curl{\begin{bmatrix}\vek C & \vek V_{m+1}\end{bmatrix}\overline{\vek G}_{m} + \sigma \begin{bmatrix}\vek Z_{\CU} & \vek Z_{m}\end{bmatrix}}^{\ast}\vek r\SfAlg{\sigma_{\ell}}$
\normalsize&\ref{line.N-inv}&2(m+k)n&L\\
Apply $\vek N^{-1}$&\ref{line.N-inv}&$\frac{2}{3}(m+k)^{3} + (m+k)^{2}$&L\\
Update approx.&\ref{line.update-x}&2(m+k)n&L\\
Update resid.&\ref{line.update-r}&2(m+k)n&L\\
\end{tabular}
\end{center}
\end{table}

\begin{table}[htb]
\caption{One-time overhead costs in Algorithm \ref{algorithm.sgcrodr} when compared to \RGMRES.}
\label{table.sgcrodr-ops-dazu-onetime}
\begin{center}
\begin{tabular}{r|c|c|c}
Operations & Alg. Line & FLOPS in $\CO(\cdot)$ & $\times$ per method execution \\
\hline
$\vek Z_{\CU} \leftarrow \vek Z_{\CU}\vek R^{-1}$ & \ref{line.Zu-update} &$k^{3}$&$1$\\
$\vek C^{\ast}\vek Z_{\CU} $&\ref{line.CTZu-ZuTZu}&$k^{2}$n&$1$\\
$\vek Z_{\CU}^{\ast}\vek Z_{\CU}$&\ref{line.CTZu-ZuTZu}&$k^{2}n$&$1$\\
Sum of $4$ $k\times k$ matrices &\ref{line.initial-N-sum}&$3k$&$L$\\
$\prn{\vek C + \sigma\vek Z_{\CU}}^{\ast}\vek r\SfAlg{\sigma_{\ell}}$&\ref{line.initN-invert-update}&$2kn$& $L$\\
Apply $\vek N^{-1} $&\ref{line.initN-invert-update}&$k^{3}+\frac{2}{3}k^{2}$&$L$\\
$\vek x\SfAlg{\sigma_{\ell}} \leftarrow \vek x\SfAlg{\sigma_{\ell}} + \vek U\vek y$&\ref{line.initx-update}&$2kn$&$L$\\
$\vek r\SfAlg{\sigma_{\ell}} \leftarrow \vek r\SfAlg{\sigma_{\ell}} - \prn{\vek C + \sigma\vek Z_{\CU}}\vek y$&\ref{line.initr-update}&$2kn$&$L$
\end{tabular}
\end{center}
\end{table}

%% file: cost_calc.tex
\begin{figure}[htb]
\hfill
\includegraphics[scale=0.23]{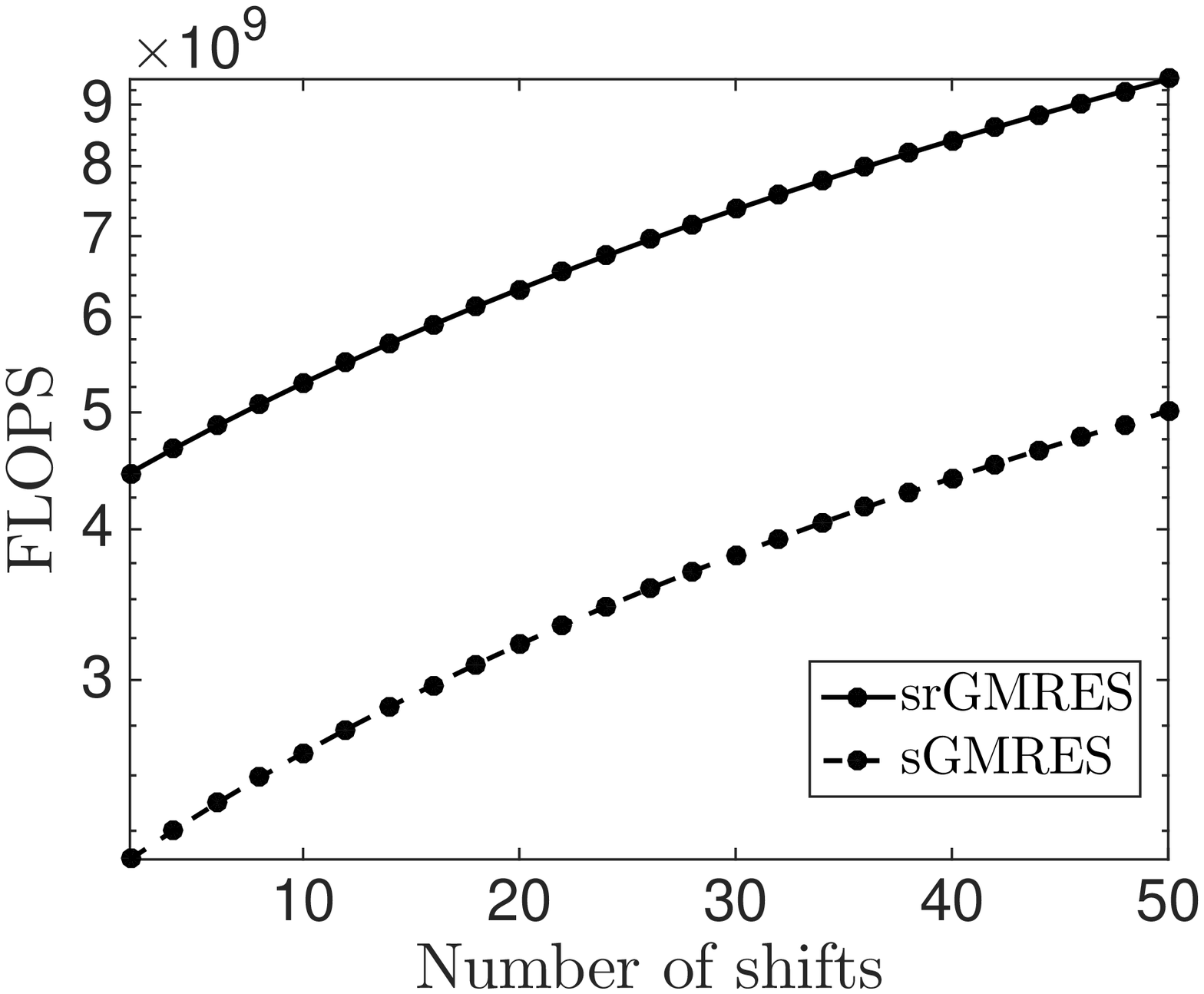}\,\,\includegraphics[scale=0.23]{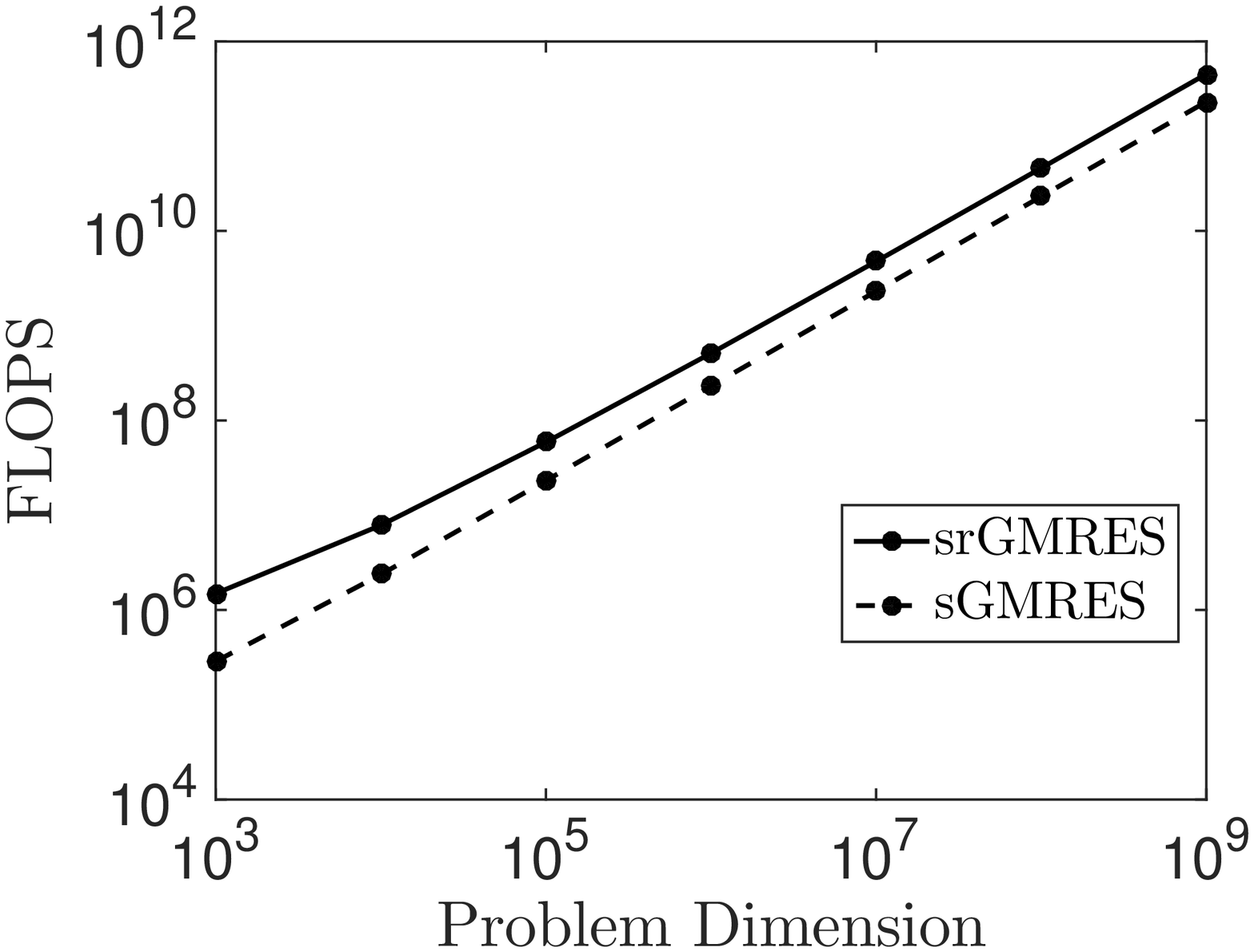}\,\,\includegraphics[scale=0.23]{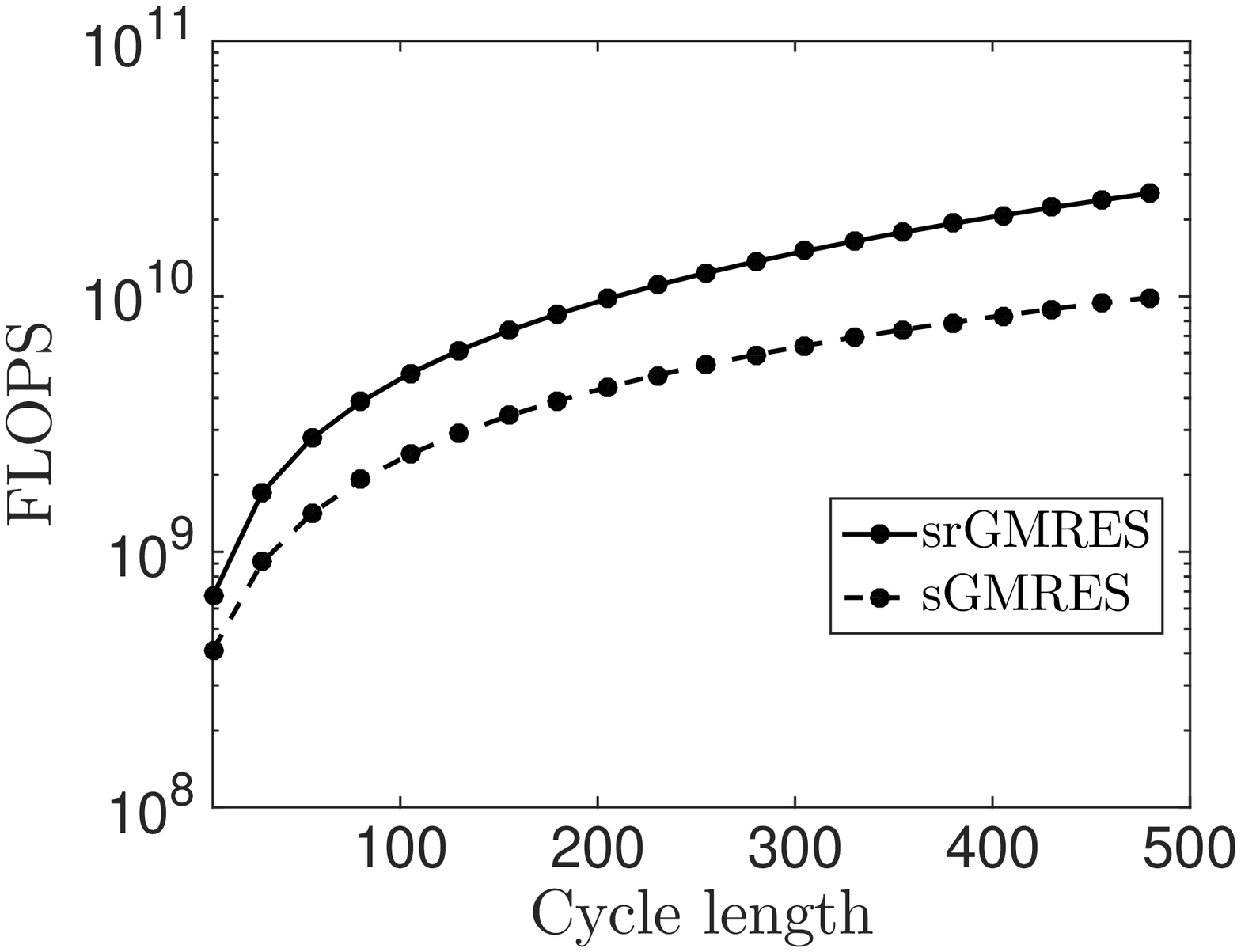}
\hfill
\begin{picture}(0,0)
\end{picture}
\caption{\label{figure.cost_calc} Estimated cost in FLOPS, respectively, for different numbers of shifts, problem dimensions, and cycle lengths with all other parameters being held constant. }
\end{figure}

%% file: qcdSmallmkTable.tex
\begin{table}[htb]
\caption{Matrix vector product counts for different pairs $(m,k)$ of restart cycle length and recycled subspace dimension for shifted \RGMRES. The matrices used in this experiment are the smaller set of QCD matrices from \cite{DH.2011}.
Experiments were performed for larger values than shown but no further improvement was observed}
\label{table.qcdSmallmkTable}
\begin{center}
\begin{tabular}{r|c|c|c|c|c|c|c|c|c}
$m\backslash k$ & $5$ & $15$ & $25$ & $35$ & $45$ & $55$ & $65$ & $75$ & $85$ \\
\hline
5 & 1566 & 1295 & 1205 & 1161 & 1146 & 1131 & 1126 & 1116 & 1111 \\
20 & 1466 & 1254 & 1182 & 1141 & 1122 & 1110 & 1107 & 1103 & 1096 \\
35 & 1418 & 1229 & 1166 & 1132 & 1113 & 1103 & 1096 & 1095 & 1091 \\
50 & 1363 & 1223 & 1158 & 1128 & 1114 & 1105 & 1099 & 1097 & 1090 \\
65 & 1344 & 1219 & 1159 & 1124 & 1109 & 1106 & 1099 & 1090 & 1086 \\
80 & 1321 & 1210 & 1153 & 1123 & 1109 & 1102 & 1098 & 1091 & 1085 \\
95 & 1321 & 1210 & 1153 & 1124 & 1108 & 1100 & 1097 & 1093 & 1084 \\
\end{tabular}
\end{center}
\end{table}

%% file: qcdLargeTimingComparison.tex
\begin{table}[htb]
\caption{Timing (in seconds) and matrix-vector product (mat-vec) comparisons between preconditionedshifted rGMRES, shifted GMRES, and sequential applications of rGMRES with cycle length $m=80$ and recycled subspace dimension $k=10$ applied to the large QCD matrices. The same preconditioner was used in all experiments.}
\label{table.qcdLargeTimingComparison}
\begin{center}
\begin{tabular}{c|c|c}
\textbf{Method} & \textbf{mat-vecs} & \textbf{time} \\
\hline
srGMRES & 3117 & 358.44 \\
sGMRES & 4003 & 322.65 \\
Seq. rGMRES & 4379 & 469.78 \\
Seq. GMRES & 5665 & 489.16 \\
\end{tabular}
\end{center}
\end{table}

%% file: qcdBigRecycleDimsPlot.tex
\begin{figure}[htb]
\hfill
\includegraphics[scale=0.35]{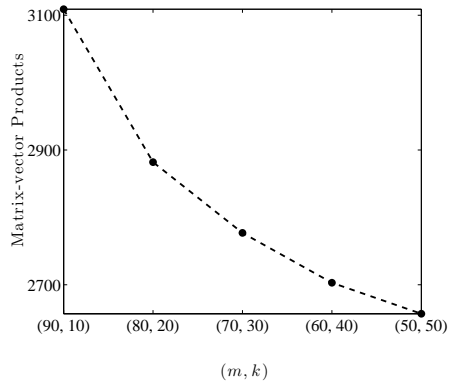}
\hfill
\begin{picture}(0,0)
\end{picture}
\caption{Matrix-vector product counts for shifted \RGMRES\ 
for the sequence of large QCD amtrices and the same shifts as in Table {\rm\ref{table.qcdLargeTimingComparison}}.
for various pairs $(m,k)$ of Krylov
subspace dimension and recycled subspace dimension such that the 
total augmented subspace Krylov subspace dimension $m+k=100$.
\label{figure.qcdBigRecycleDimsPlot}
}
\end{figure}

%% file: qcdLargeShiftSizeMatvecCounts.tex
\begin{table}[htb]
\caption{Comparison of 3 methods for different shifts sizes. In each experiment, two systems were solved, the base system and one shifted system with the shift shown in the table column header.}
\label{table.qcdLargeShiftSizeMatvecCounts}
\begin{center}
\begin{tabular}{c|c|c|c|c|c|c|c}
Method\textbackslash $\norm{\sigma}$ & $10^{-3}$ & $10^{-2}$ & $10^{-1}$ & $10^{0}$ & $10^{1}$ & $10^{2}$ & $10^{3}$ \\
\hline
Sh. GMRES Alg. \ref{algorithm.sgmres} & 1330 & 1405 & 1294 & 967 & 1067 & 1265 & 1306 \\
Sh. rGMRES Alg. \ref{algorithm.sgcrodr} & 980 & 1039 & 1017 & 804 & 908 & 1105 & 1144 \\
Seq. rGMRES & 1183 & 1170 & 1077 & 812 & 914 & 1111 & 1152 \\
\end{tabular}
\end{center}
\end{table}

%% file: qcdLargeMarginalCost.tex
\begin{figure}[htb]
\hfill
\includegraphics[scale=0.35]{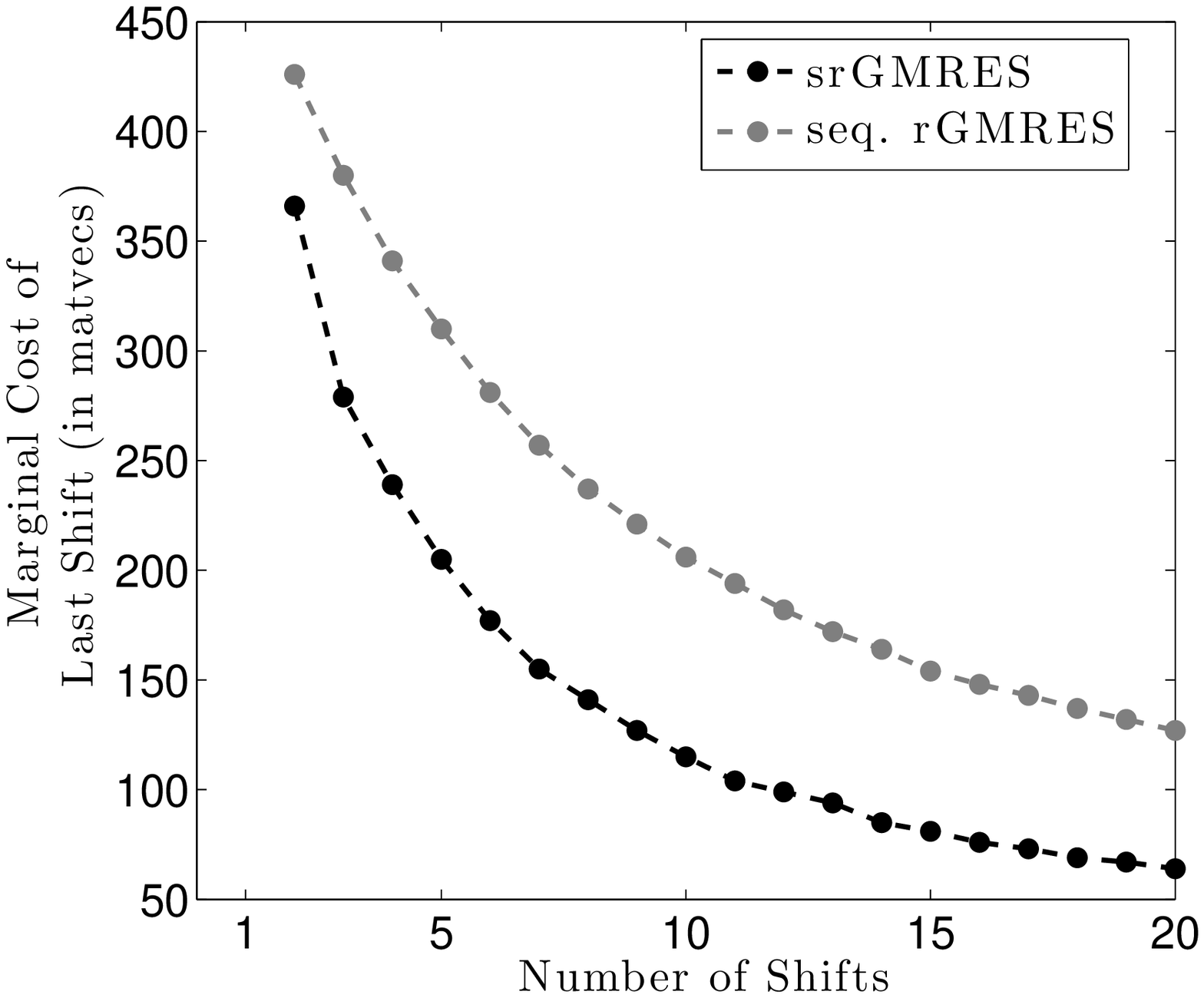}
\includegraphics[scale=0.35]{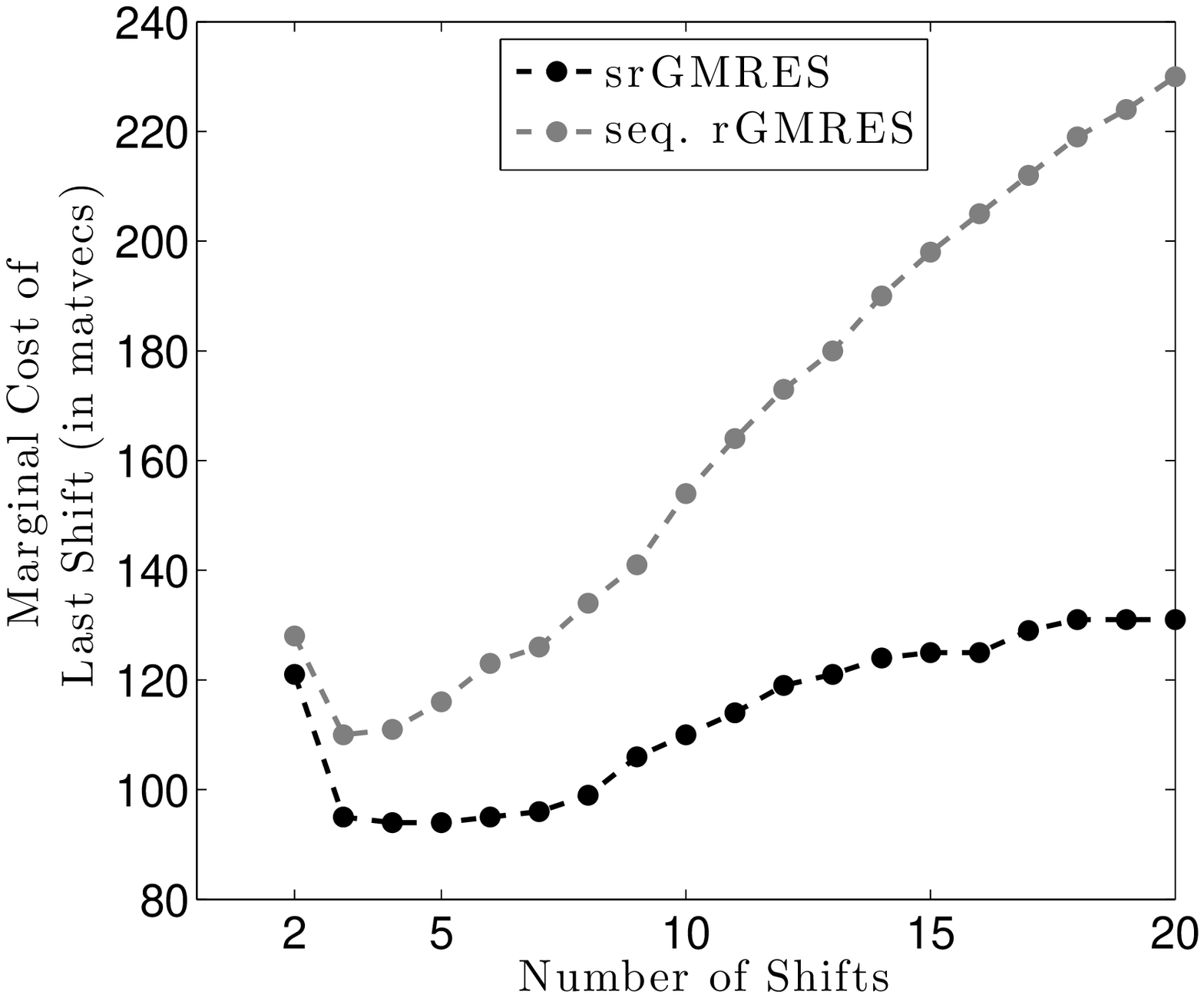}
\hfill
\begin{picture}(0,0)
\end{picture}
\caption{Comparison of the marginal cost of solving each addition shifted system.
For the left-hand figure, the shifts were evenly space points
from the interval $\brac{0,1}$, and in the right-hand figure, the shifts were evenly spaced points
from the larger interval $\brac{1,10}$
\label{figure.qcdLargeMarginalCost}
}
\end{figure}

%% file: shiftMagnitudeConvergenceComparison.tex
\begin{figure}[htb]
\hfill
\includegraphics[scale=0.60]{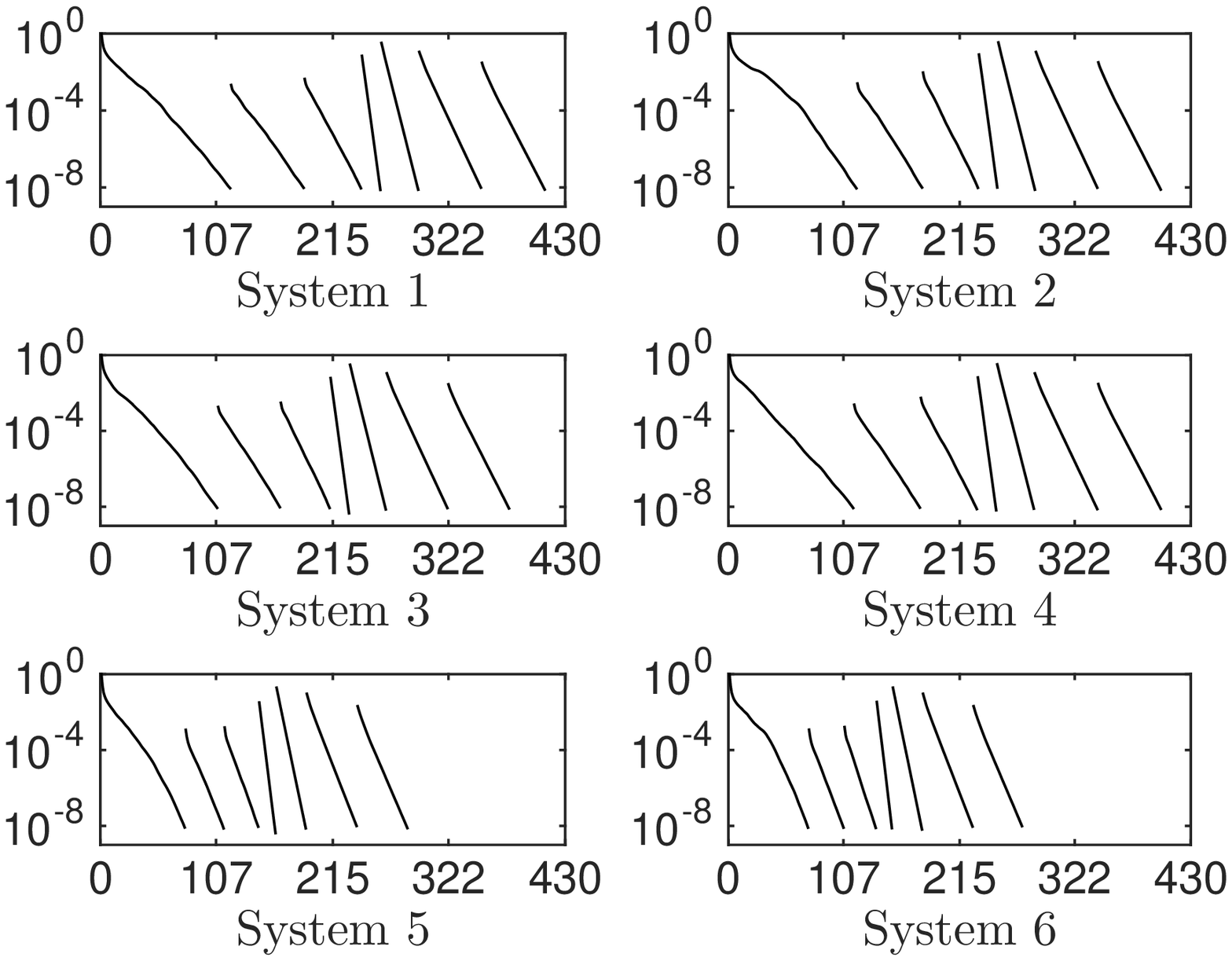}
\hfill
\begin{picture}(0,0)
\end{picture}
\caption{For the large QCD matrices and $(m,k)=(100,5)$, an illustration of the
amount of residual improvement for different magnitude shifts,
$\sigma\in\curl{10^{-3},10^{-2},\ldots, 10^{3}}$.  In each subplot, we display
the residual curves sequentially to reflect that the algorithm is called for each shifted
linear system in sequence.  The order in which the systems were solved is the same as the
order of the listed shifts.  
\label{figure.shiftMagnitudeConvergenceComparison}
}
\end{figure}